\documentclass[11pt]{amsart}
\usepackage{amsfonts}
\usepackage{amssymb, amstext, amscd, amsmath}
\usepackage{amsthm}
\usepackage{times}
\usepackage{fullpage}
\usepackage{txfonts}
\usepackage{indentfirst}
\usepackage[all]{xy}
\usepackage{subfigure}
\usepackage{graphicx}
\usepackage{tikz}
\usepackage{pifont}
\usepackage{bbm}
\usepackage{amsfonts}
\usepackage{mathrsfs}
\usepackage{lscape}
\usepackage{setspace} 
\usepackage{graphicx}      
\usepackage{subfigure}     
\usepackage{color}              
\usepackage{indentfirst}        
\usepackage{url}
\usepackage{fancyhdr}
\usepackage[nocompress]{cite}
\newtheoremstyle{noparens}%
  {}{}%
  {\itshape}{}%
  {\bfseries}{.}%
  { }%
  {\thmname{#1}\thmnumber{ #2}\mdseries\thmnote{ #3}}
\theoremstyle{noparens}
\allowdisplaybreaks

\newtheorem{theorem}{Theorem}[section]

\newtheorem{proposition}{Proposition}[section]
\newtheorem{lemma}{Lemma}[section]
\newtheorem{remark}{Remark}[section]
\newtheorem{definition}{Definition}[section]
\newtheorem{example}{Example}[section]

\numberwithin{equation}{section}
\numberwithin{equation}{section}

\begin{document}
\baselineskip 18pt

\title[The $F$-polynomial invariant for knotoids]{The $F$-polynomial invariant for knotoids}

\author{Yi Feng}
\address{School of Mathematical Sciences, Dalian University of Technology, Dalian 116024, P. R. China}
\email{fengyi@mail.dlut.edu.cn}

\author{Fengling Li$^*$}
\address{School of Mathematical Sciences, Dalian University of Technology, Dalian 116024, P. R. China}
\email{fenglingli@dlut.edu.cn}
\thanks{$*$Supported in part by a grant (No. 12071051) of NSFC}


\begin{abstract}
As a generalization of the classical knots, knotoids deal with the open ended knot diagrams in a surface. In recent years, many polynomial invariants for knotoids have appeared, such as the bracket polynomial, the index polynomial and the $n$th  polynomial, etc. In this paper, we introduce a new polynomial invariant $F$-polynomial for knotoids and discuss some properties of the $F$-polynomial. Then we construct a family of knotoid diagrams which can be distinguished from each other by the $F$-polynomial but cannnot be distinguished by the index polynomial and the $n$th  polynomial.
\end{abstract}
\maketitle

\section{Introduction}\label{introduction }
 Knotoid theory was first proposed by Turaev \cite{16} in 2012.
 A knotoid diagram in a surface $\Sigma$ is a generic immersion from the unit interval $[0,1]$ to $\Sigma$ with finitely many transversal double points endowed with over and under crossings information, the images of 0 and 1 are the endpoints of the knotoid diagram.  A knotoid is the equivalence class of knotoid diagrams up to istopies of $\Sigma$ and classical Reidemeister moves of $R_1$, $R_2$ and $R_3$ away from the endpoints. In this paper, we only consider the knotoids in $\Sigma=S^2$ which are called spherical knotoids. For a given knotoid diagram, we  can also find the Gauss diagram of it.

A (classical) knotoid diagram in $S^2$ is an open-ended knot diagram. It is the generalization of
the long knot diagram, while the long knot diagram only allows its endpoints to be  in the same region, the knotoid diagram allows its endpoints to be in different regions of the diagram. The theory of knotoid diagrams suggests a new diagrammatic approach to knots. Turaev introduced the method of converting knotoids into knots in \cite{16} which was later called underpass/overpass closure map by G\"{u}g\"{u}mc\"{u} and Kauffman in \cite{7}, and they proved that it is surjective rather than injective. Every knotoid diagram $D$ in $S^2$ determines a knot in $S^3$, it can be obtained by connecting the endpoints of $D$  with an embedded arc $a\subset S^2$ which meets $D$ transversely at a finite set of points distinct from the crossings of $D$, and  passes everywhere under/over $D$. Conversely, given a knot in $S^3$ we also have  a knotoid diagram corresponding to it. In fact, we only need to cut out an underpassing/overpassing strand. The strand may contain no crossing, 1 crossing, or even more. See \cite{16} for more details.

 Knotoid theory has been getting a lot of attention since it was proposed. Knotoid invariant is an important topic in knotoid theory. In \cite{16}, Turaev  introduced knotoids, gave the relationship between knotoids and knots and extended many knot invariants to knotoids, such as the knot group, the bracket polynomial. In \cite{7}, G\"{u}g\"{u}mc\"{u} and Kauffman introduced the theory of virtual knotoids, defined some closure maps from knotoids to classical knots and virtual knots and constructed the odd writhe, the parity bracket polynomial, the affine index polynomial and the arrow polynomial for knotoids. In \cite{13}, Miyazawa  used  the concept of the pole diagram to define the extended bracket polynomial which is an enhancement of the bracket polynomial for knotoids. And  Kim et al.\cite{12}  introduced the index polynomial and the $n$th polynomial for  (virtual) knotoids. In addition, some algebraic invariants have also been given in \cite{8,9,10}, such as biquandle coloring invariants, biquandle bracket matrix, quantum invariant and so on.
  In \cite{14}, Manousos et al.  defined finite type invariants for knotoids via the Vassiliev skein relation, and they showed that there are non-trivial type-1 invariants for spherical knotoids.

In the past few years, knotoids have been used to study the structure of open protein chains, the biological purpose of the knotoids in proteins is  an interesting topic. In \cite{6}, G\"{u}g\"{u}mc\"{u} et al.  studied bonded knotoids $(K,b)$ where $K$ is a knotoid and $b$ is a finite collection of arcs connecting two distinct interior points of $K$, constructed topological invariants of bonded knotoids and gave an application to indicate that these invariants can be used to analyze the topological structure of proteins. In \cite{5},  Goundaroulis et al.  used planar knotoids and bonded knotoids to provide a more detailed overview of the topology of open protein chains and set up a powerful topology model that can distinguish proteins with disulphide. In \cite{4}, Goundaroulis et al.  used the concept of knotoids to analyze and reveal that the knotoid approach can detect protein regions that were classified earlier as knotted. In \cite{2}, Barbensi and Goundaroulis determined an optimal size of sampled projections set using the  $f$-distance $d_f(k_1,k_2)$ which is the minimal number of forbidden moves to deform $k_1$ to $k_2$, probed the depth of knotted proteins with trefoil as main knot type and proved the knotoid spectrum of a knotted protein contains important information on the geometry and topology of the protein itself.

Recall the definitions of the index polynomial and the $n$th polynomial which are invariants for (virtual) knotoids in \cite{12}.
The index polynomial is defined by
$$
F_{D}(t)=\sum_{c\in C(G(D))}{\rm sign}(c)(t^{i(c)}-1),
$$
where $G(D)$ is the  Gauss diagram of $D$, $C(G(D))$ is the set of chords of $G(D)$ and $i(c)$ is the intersection index of $c$.
For each nonnegative integer $n$, the $n$th  polynomial is defined by
$$
Z_{D}^{n}(t)=\sum_{c\in C_{n}(G(D))}{\rm sign}(c)(t^{d_n(c)}-1),
$$
where $C_{n}(G(D))= \left\{
c\in C(G(D))|i(c)=kn~ {\rm for~ some~integer~}~k
\right\}$, $d_n(c)$ is the sum of signs of the endpoints on the lead side of $c$ whose chords in $C_{n}(G(D))\backslash
\left\{
 c
 \right\}$.
In particular, $Z_{D}^{1}(t)=F_D(t).$

 We observe that there exist knotoid diagrams that the index polynomial and the $n$th polynomial fail to distinguish. In this paper, we give a new polynomial invariant $F_D(u,v)$ for knotoids  which is called $F$-$polynomial$.
Let $D$ be a knotoid diagram, $G(D)$  the Gauss diagram of $D$ and $C(G(D))$ the set of chords of $G(D)$. The $F$-polynomial is defined as follows:
$$
F_D(u,v)=\sum_{c\in C(G(D))}{\rm sign}(c)(u^{g_c(v)}-1).
$$
Then we discuss some properties of it and give a family of knotoid diagrams that can be distinguished by $F$-polynomial.

 The paper is organized as follows. In Section 2, we recall the relevant definitions and conclusions for knotoids. In Section 3, we define a polynomial invariant $F_D(u,v)$ for knotoids  by using the Gauss diagram of knotoid diagram $D$, give some properties of it and generalize it to a class of polynomial invariants for knotoids. In Section 4, we give a family of knotoid diagrams which can be distinguished from each other by $F$-polynomial but cannot be distinguished by the index polynomial and the $n$th  polynomial.

\section{Preliminaries}\label{Preliminary}
In this section, we give some definitions and conclusions about  knotoids \cite{7,12,16}.
\begin{definition}
{\rm A $knotoid ~diagram$ is a generic immersion from the unit interval [0,1]
to $S^2$ with finitely many transversal double points endowed with over
and under crossings information. These points are called $classical ~crossings$ of the knotoid diagram. The $endpoints$ of the knotoid diagram are the images of 0 and 1 which are different from the double points and from each other, and they are called $the ~tail$ and $the ~head$ of the knotoid diagram, respectively.
The orientation of a knotoid diagram is from the tail to the head.}

{\rm A $trivial~knotoid ~diagram$ is an embedding from the unit interval $[0,1]$
to $S^2$.}
\end{definition}

\begin{figure}
  \centering
  \begin{tikzpicture}[scale=1]
\draw [thick](0,1)--(0,3);
\draw[thick,<->](0.25,2)--(0.75,2);
\draw[thick](1,3)--(3,1);
\draw[thick](1,1)--(1.9,1.9);
\draw[thick](2.1,2.1)--(3,3);
\draw[thick](3,1) to [out=75,in=285](3,3);

\node at (1,0){$R_1$};
\end{tikzpicture}
\qquad
\qquad
\begin{tikzpicture}[scale=1]
\draw [thick](0,1)--(0,3);
\draw [thick](0.5,1)--(0.5,3);
\draw[thick,<->](0.75,2)--(1.25,2);
\draw[thick](1.5,1)--(2.5,2)--(1.5,3);
\draw[thick](2.5,1)--(2.1,1.4);
\draw[thick](1.9,1.6)--(1.5,2)--(1.9,2.4);
\draw[thick](2.1,2.6)--(2.5,3);
\node at (1,0){$R_2$};

\end{tikzpicture}
\qquad
\qquad
\begin{tikzpicture}[scale=1]
\draw [thick](0,1)--(2,3);
\draw [thick](0,3)--(0.9,2.1);
\draw[thick](1.1,1.9)--(2,1);
\draw[thick](0,2.5)--(0.4,2.5);
\draw[thick](0.6,2.5)--(1.4,2.5);
\draw[thick](1.6,2.5)--(2,2.5);
\draw[thick,<->](2.25,2)--(2.75,2);
\draw[thick](3,1)--(5,3);
\draw[thick](3,1.5)--(3.4,1.5);
\draw[thick](3.6,1.5)--(4.4,1.5);
\draw[thick](4.6,1.5)--(5,1.5);
\draw[thick](3,3)--(3.9,2.1);
\draw[thick](4.1,1.9)--(5,1);
\node at (2.5,0){$R_3$};
\end{tikzpicture}
\caption{Classical Reidemeister moves.}
\label{Reidemeister move}
  \end{figure}
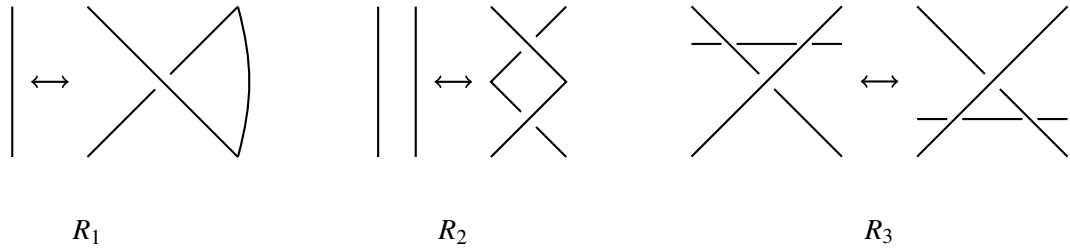

Just like knots, for each crossing of knotoid diagrams, we have two types of crossing signs as shown in Figure~\ref{crossing signs}.
\begin{figure}
  \centering
 \subfigure [sign($c$)=$+1$]{
  \begin{tikzpicture}[scale=1]
\draw[thick,->](0,0)--(2,2);
\draw[thick](2,0)--(1.1,0.9);
\draw[thick,->](0.9,1.1)--(0,2);
\node [below] at (1,0.9){$c$};
\end{tikzpicture}}
\qquad
\qquad
\qquad
\qquad
\qquad
\qquad
\subfigure [sign($c$)=$-1$]{
\begin{tikzpicture}[scale=1]
\draw[thick,<-](0,2)--(2,0);
\draw[thick](0,0)--(0.9,0.9);
\draw[thick,->](1.1,1.1)--(2,2);
\node[below] at (1,0.9){$c$};
\end{tikzpicture}}
  \caption{The sign of crossing $c$.}
  \label{crossing signs}
\end{figure}
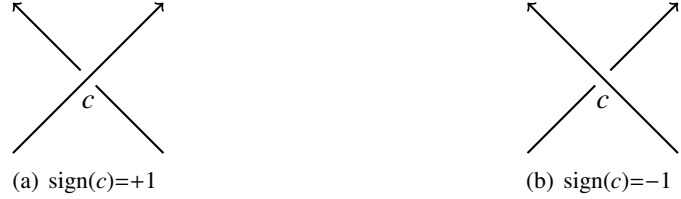
\begin{definition}
{\rm Two knotoid diagrams are said to be $equivalent$ if one of them can be deformed to
the other by isotopies of $S^2$ and a finite sequence of classical Reidemeister moves $R_1$, $R_2$ and $R_3$ away from the endpoints as shown in Figure~\ref{Reidemeister move}.}
\end{definition}

   A $knotoid$ is the equivalence class of knotoid diagrams by isotopies of $S^2$ and  classical Reidemeister moves taking place in a local disk free of the endpoints.

\begin{figure}
  \begin{tikzpicture}[scale=1]
\draw [thick](0,1)--(0,3);
\draw[thick,<->](0.25,2)--(0.75,2);
\draw[thick](1,3)--(3,1);
\draw[thick](1,1)--(1.9,1.9);
\draw[thick](2.1,2.1)--(3,3);
\draw[thick](3,1) to [out=75,in=285](3,3);

\node at (1,0){$R_1$};
\end{tikzpicture}
\qquad
\qquad
\begin{tikzpicture}[scale=1]
\draw [thick](0,1)--(0,3);
\draw [thick](0.5,1)--(0.5,3);
\draw[thick,<->](0.75,2)--(1.25,2);
\draw[thick](1.5,1)--(2.5,2)--(1.5,3);
\draw[thick](2.5,1)--(2.1,1.4);
\draw[thick](1.9,1.6)--(1.5,2)--(1.9,2.4);
\draw[thick](2.1,2.6)--(2.5,3);
\node at (1,0){$R_2$};

\end{tikzpicture}
\qquad
\qquad
\begin{tikzpicture}[scale=1]
\draw [thick](0,1)--(2,3);
\draw [thick](0,3)--(0.9,2.1);
\draw[thick](1.1,1.9)--(2,1);
\draw[thick](0,2.5)--(0.4,2.5);
\draw[thick](0.6,2.5)--(1.4,2.5);
\draw[thick](1.6,2.5)--(2,2.5);
\draw[thick,<->](2.25,2)--(2.75,2);
\draw[thick](3,1)--(5,3);
\draw[thick](3,1.5)--(3.4,1.5);
\draw[thick](3.6,1.5)--(4.4,1.5);
\draw[thick](4.6,1.5)--(5,1.5);
\draw[thick](3,3)--(3.9,2.1);
\draw[thick](4.1,1.9)--(5,1);
\node at (2.5,0){$R_3$};

\end{tikzpicture}
\qquad
\qquad
\\

\begin{tikzpicture}[scale=1]
\draw [thick](0,1)--(0,3);
\draw[thick,<->](0.25,2)--(0.75,2);
\draw[thick](1,3)--(3,1);
\draw[thick](1,1)--(3,3);
\draw[thick](3,1) to [out=75,in=285](3,3);
\draw[thick](2,2) circle[radius=0.1];
\node at (1,0){$VR_1$};
\end{tikzpicture}
\qquad
\qquad
\begin{tikzpicture}[scale=1]
\draw [thick](0,1)--(0,3);
\draw [thick](0.5,1)--(0.5,3);
\draw[thick,<->](0.75,2)--(1.25,2);
\draw[thick](1.5,1)--(2.5,2)--(1.5,3);
\draw[thick](2.5,1)--(1.5,2)--(2.5,3);
\draw[thick](2.05,2.55)--(2.5,3);
\draw[thick](2,1.5) circle[radius=0.1];
\draw[thick](2,2.5) circle[radius=0.1];
\node at (1,0){$VR_2$};
\end{tikzpicture}
\qquad
\qquad
\begin{tikzpicture}[scale=1]
\draw [thick](0,1)--(2,3);
\draw [thick](0,3)--(2,1);
\draw[thick](0,2.5)--(2,2.5);
\draw[thick](0.5,2.5) circle[radius=0.1];
\draw[thick](1.5,2.5) circle[radius=0.1];
\draw[thick](1,2) circle[radius=0.1];
\draw[thick,<->](2.25,2)--(2.75,2);
\draw[thick](3,1)--(5,3);
\draw[thick](3,1.5)--(5,1.5);
\draw[thick](3,3)--(5,1);
\draw[thick](3.5,1.5) circle[radius=0.1];
\draw[thick](4.5,1.5) circle[radius=0.1];
\draw[thick](4,2) circle[radius=0.1];
\node at (2.5,0){$VR_3$};
\end{tikzpicture}
\qquad
\qquad
\\
\begin{tikzpicture}[scale=1]
\draw [thick](0,1)--(2,3);
\draw [thick](0,3)--(0.9,2.1);
\draw[thick](0,2.5)--(2,2.5);
\draw[thick](1.1,1.9)--(2,1);
\draw[thick](0.5,2.5) circle[radius=0.1];
\draw[thick](1.5,2.5) circle[radius=0.1];
\draw[thick,<->](2.25,2)--(2.75,2);
\draw[thick](3,1)--(5,3);
\draw[thick](3,1.5)--(5,1.5);
\draw[thick](3,3)--(3.9,2.1);
\draw[thick](4.1,1.9)--(5,1);
\draw[thick](3.5,1.5) circle[radius=0.1];
\draw[thick](4.5,1.5) circle[radius=0.1];
\node at (2.5,0){$VR_4$};
\end{tikzpicture}
\qquad
\qquad
\qquad
\begin{tikzpicture}[scale=1.2]
  \draw[thick](1.25,0.8)--(2.75,0.8);
  \fill (2.75,0.8) circle [radius=0.075];
  \draw[thick](4,1.6) arc (90:270:0.8);
  \draw[thick,<->] (4.5,0.8)--(5,0.8);
  \draw[thick](5.5,0.8)--(7,0.8);
  \fill (7,0.8) circle [radius=0.075];
  \draw[thick](7,1.6) arc (90:270:0.8);
  \draw[thick](6.2,0.8)circle[radius=0.1];
  \node at (4.75,-0.8){$\Omega_v$};
\end{tikzpicture}
\caption{Generalized Reidemeister moves and {$\Omega_v$}-move.}
\label{moves}
  \end{figure}
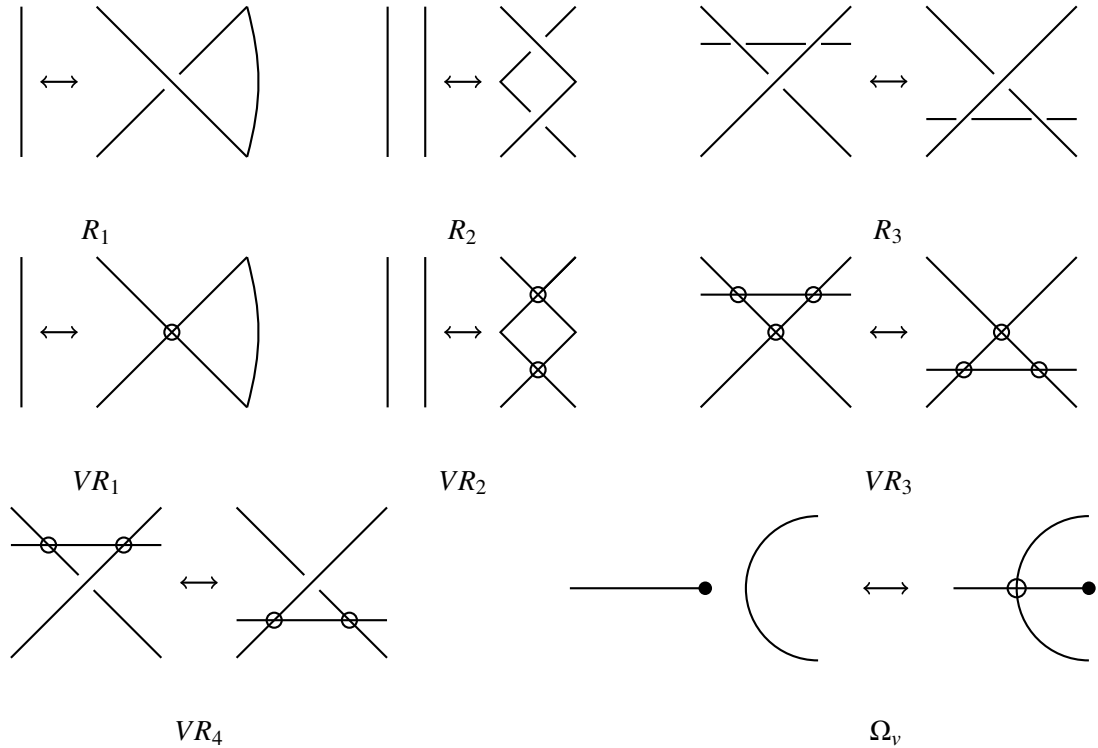

    Depending on the location of the endpoints, there are two types of knotoids. One is $knot$-$type~knotoid$ which has a knotoid diagram that the endpoints lie in the same region, the other is $proper ~knotoid $ which has all knotoid diagrams that the endpoints are not in the same region.
  \begin{remark}   There are two types of moves  $\phi_+$ and $\phi_-$ as shown in Figure~\ref{p} that are not allowed. They are called the $forbidden ~moves$. In fact, if they are allowed, then any knotoid can be deformed into trivial knotoid.
  \end{remark}
  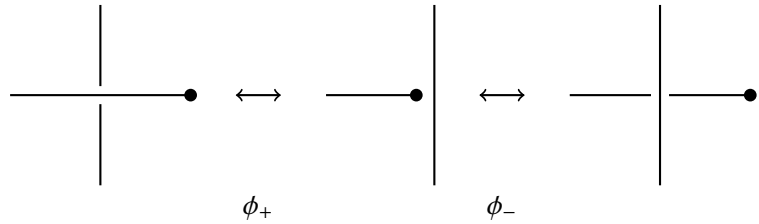
\begin{figure}
    \centering
\begin{tikzpicture}[scale=1.2]
\draw[thick](-2.5,0)--(-0.5,0);
\draw[thick](-1.5,-1)--(-1.5,-0.1);
\draw[thick](-1.5,0.1)--(-1.5,1);
\draw [thick,<->](0,0)--(0.5,0);
\fill (-0.5,0) circle [radius=2pt];
\draw[thick](1,0)--(2,0);
\draw[thick](2.2,-1)--(2.2,1);
\draw [thick,<->](2.7,0)--(3.2,0);
\fill (2,0) circle [radius=2pt];
\draw[thick](3.7,0)--(4.6,0);
\draw[thick](4.8,0)--(5.7,0);
\draw[thick](4.7,-1)--(4.7,1);
\fill (5.7,0) circle [radius=2pt];
\node [below]at(0.25,-1){$\phi_+$};
\node [below]at(2.95,-1){$\phi_-$};

\end{tikzpicture}
    \caption{Forbidden moves.}\label{p}
  \end{figure}
\begin{definition}
{\rm A $virtual~knotoid~diagram$ is the knotoid diagram that may contain two types of crossings, one is classical crossing and the other is the crossing with enclosed circle without over/under information. It is called virtual crossing  as shown in Figure~\ref{crossing type}~\subref{5}.}
\end{definition}
\begin{figure}
\centering
 \subfigure [classical crossing]{
\begin{tikzpicture}[scale=1.3]
\draw[thick](0,2)--(2,0);
\draw[thick](0,0)--(0.9,0.9);
\draw[thick](1.1,1.1)--(2,2);
\end{tikzpicture}}
\qquad
\qquad
\qquad
 \subfigure [virtual crossing]{
\begin{tikzpicture}[scale=1.3]
\draw[thick](0,2)--(2,0);
\draw[thick](0,0)--(2,2);
\draw[thick](1,1) circle[radius=0.1];
\label{5}
\end{tikzpicture}}
\caption{Two types of crossings.}
\label{crossing type}
\end{figure}
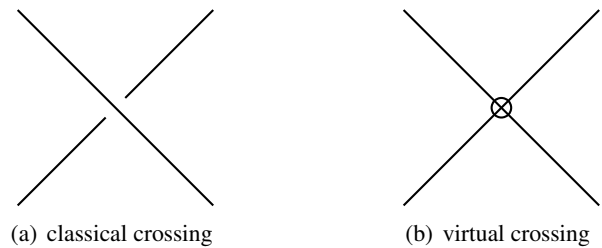
\begin{definition}
{\rm Two virtual knotoid diagrams are said to be $equivalent$ if one of them can be deformed to
the other by a finite sequence of Generalized Reidemeister moves and {$\Omega_v$}-moves away from the endpoints as shown in Figure~\ref{moves}.}
\end{definition}

  A $virtual~knotoid$ is the equivalence class of virtual knotoid diagrams up to Generalized Reidemeister moves and {$\Omega_v$}-moves taking place in a local disk free of the endpoints.

\begin{figure} 
\centering
\begin{tikzpicture}[scale=0.7]
    \draw[thick](0,0) arc (30:330:2);
    \draw[thick,dashed] (0,-1.99) arc (330:390:2);
    \fill (0,0) circle (.1);
    \node at (0.5,0) [right]{the tail};
    \fill (0,-2) circle (.1);
    \node at (0.5,-2) [right]{the head};
    \draw [thick,->](-1.7,1) -- (-1.7,-3);
    \node [above] at(-1.5,1){sign($c$)};
    \node [below]at(-1.5,-3){$-$sign($c$)};
    \node[right] at (-1.5,-1){$c$};
    \node [below]at(-1.5,-4){$G(D)$};
\end{tikzpicture}
\caption{Gauss diagram of $D$.}
\label{Gauss diagram}
\end{figure}
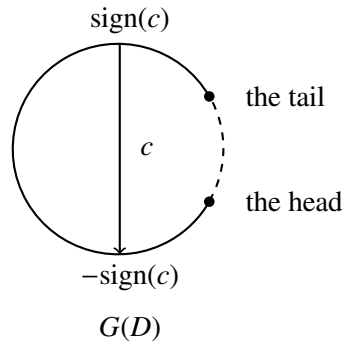

 \begin{definition}
 {\rm Let $D$ be a knotoid diagram with $m$ classical crossings. A $Gauss ~diagram$ of $D$ is a counter-clockwise oriented arc $C$ where the preimages of the overpassing strand and the underpassing strand of each crossing are connected by a chord. Chords are directed from the preimages of the overpassing strands to the preimages of the underpassing strands. Since there is a bijection between the classical crossings of $D$ and the chords of the Gauss diagram of $D$, we still write $c$ for the chord corresponding to each crossing $c$. We give sign($c$) to the preimage of the overpassing strand of $c$ and $-$sign($c$) to the preimage of the underpassing strand of $c$. The Gauss diagram of $D$ is denoted by $G(D)$. The tail(head) of $D$ is the initial(terminal) point of $G(D)$ which is called $the~ tail(head)$ of $G(D)$, respectively, as shown in Figure~\ref{Gauss diagram}. Note that there is no chord on the dotted arc between the head and the tail.}
 \end{definition}

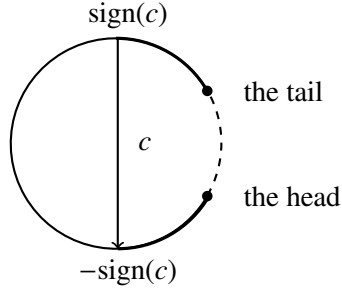
\begin{figure}
  \centering
   \begin{tikzpicture}[scale=0.7]
    \draw[thick](0,0) arc (30:330:2);
    \draw[thick,dashed] (0,-1.99) arc (330:390:2);

    \node at (0.5,0) [right]{the tail};

    \node at (0.5,-2) [right]{the head};
    \draw [thick,->](-1.7,1) -- (-1.7,-3);
    \draw [very thick] (0,0) arc (30:90:2);
     \draw [very thick] (-1.7,-3) arc (270:330:2);

      \fill (0,0) circle (.1);
      \fill (0,-2) circle (.1);
    \node [above] at(-1.5,1){sign($c$)};
    \node [below]at(-1.5,-3){$-$sign($c$)};
    \node[right] at (-1.5,-1){$c$};
\end{tikzpicture}
  \caption{The lead side of $c$.}
  \label{lead side}
  \end{figure}

 Let $c$ be a chord of $G(D)$. The counter-clockwise oriented arc $C$ is divided into two sides by the chord $c$, and we call the side containing the head and the tail as $the~ lead~ side~ of~ c$ as shown in Figure~\ref{lead side}.
 In \cite{12}, the intersection index of $c$ was denoted by $i(c)$, which represents the sum of signs of the endpoints on the lead side of  $c$ except the signs of the endpoints of $c$. $i(c)$ was used to define the index polynomial for knotoids and a family of polynomial invariants the $n$th  polynomial for virtual knotoids.

\section{Knotoid invariants}\label{Knotoid invariant}

Let $D$ be a knotoid diagram, $G(D)$  the Gauss diagram of $D$ and $c$ the chord of $G(D)$. Denote the set of chords that pass through $c$ and point to the lead side of  $c$ by $l(c)=\left\{l_1(c), l_2(c),{\cdots},l_m(c) \right\}$ and the set of chords that point from the lead side of  $c$ to the other side by $r(c)=\left\{r_1(c), r_2(c),{\cdots},r_n(c) \right\}$ as shown in Figure~\ref{the type crossing c}. For each chord $c$, take a map $\phi_c$ from $\mathbb{Z}$ to $\mathbb{Z}_{|i(c)|}$, which maps the intersection indices of the chords that pass through  $c$ to $\mathbb{Z}_{|i(c)|}$.
\begin{figure}
  \centering
\begin{tikzpicture}[scale=0.7]
    \draw[thick](0,0) arc (30:330:2);
    \draw[thick,dashed] (0,-1.99) arc (330:390:2);
    \fill (0,0) circle (.1);
    \node at (0.5,0) [right]{the tail};
    \fill (0,-2) circle (.1);
    \node at (0.5,-2) [right]{the head};
    \draw [thick,->](-1.7,1) -- (-1.7,-3);
    \node [above] at(-1.5,1){sign($c$)};
    \node [below]at(-1.5,-3){$-$sign($c$)};
    \node[right] at (-1.5,-1){$c$};
    \draw [thick,->](-0.4,0.5) -- (-3.03,0.5);
    \draw [thick,<-](-0.44,-2.5) -- (-3.05,-2.5);
     \node [below]at(-1,0.5){$r(c)$};
     \node [above] at(-1,-2.5){$l(c)$};
\end{tikzpicture}
  \caption{The types of the chords that pass through $c$.}\label{the type crossing c}
\end{figure}
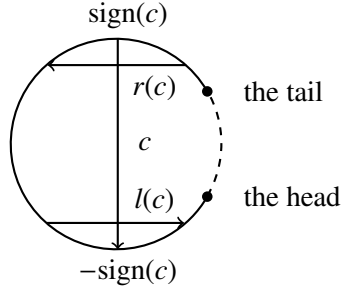
\begin{definition}
{\rm Let $D$ be a knotoid diagram and $G(D)$  the  Gauss diagram of $D$, for each chord $c$ in   $G(D)$, the $index ~function$ of $c$ is  defined by}
\begin{equation}
g_c(v)=\sum^n_{i=1}{\rm sign}(r_i(c))v^{\phi_c(i(r_i(c)))}-\sum^m_{j=1}{\rm sign}(l_j(c))v^{\phi_c(-i(l_j(c)))}.
\end{equation}
\end{definition}
In general, $g_c(v)$ takes values in  $\mathbb{Z}[v,v^{-1}]/(v^{|i(c)|}-1)$.

Then we define the $F$-polynomial as follows.
\begin{definition}
{\rm Let $D$ be a knotoid diagram, $G(D)$  the  Gauss diagram of $D$ and $C(G(D))$ the set of chords of $G(D)$. Then the $F$-$polynomial$ $F_D(u,v)$ is defined by}
\begin{equation}
F_D(u,v)=\sum_{c\in C(G(D))}{\rm sign}(c)(u^{g_c(v)}-1).
\end{equation}
\end{definition}

\subsection{Invariance of $F$-polynomial $F_D(u,v)$  }
 In this subsection, we will show that $F$-polynomial is an invariant for knotoids.
\begin{theorem}
\label{invariance theorem}

The $F$-polynomial $F_D(u,v)$ is an invariant for knotoids.
\end{theorem}

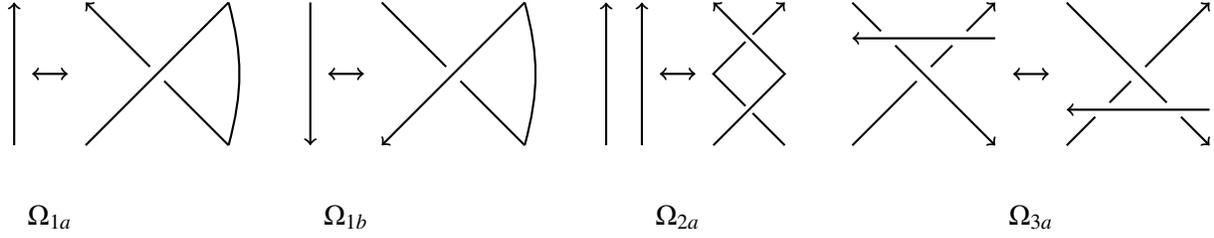
\begin{figure}
\centering
\begin{tikzpicture}[scale=0.95]
\draw[thick,->](0,1)--(0,3);
\draw[thick,<->](0.25,2)--(0.75,2);
\draw[thick,<-](1,3)--(1.9,2.1);
\draw[thick](2.1,1.9)--(3,1);
\draw[thick](1,1)--(3,3);
\draw[thick](3,1) to [out=75,in=285](3,3);
\node  at (0.5,0){$\Omega_{1a}$};
\end{tikzpicture}
\qquad
\begin{tikzpicture}[scale=0.95]
\draw[thick,<-](0,1)--(0,3);
\draw[thick,<->](0.25,2)--(0.75,2);
\node at(0.5,0){$\Omega_{1b}$};
\draw[thick](1,3)--(1.9,2.1);
\draw[thick](2.1,1.9)--(3,1);
\draw[thick,<-](1,1)--(3,3);
\draw[thick](3,1) to [out=75,in=285](3,3);
\end{tikzpicture}
\qquad
\begin{tikzpicture}[scale=0.95]
\draw[thick,->](0,1)--(0,3);
\draw[thick,->](0.5,1)--(0.5,3);
\draw[thick,<->](0.75,2)--(1.25,2);
\node at(1,0){$\Omega_{2a}$};
\draw[thick,->](1.5,1)--(2.5,2)--(1.5,3);
\draw[thick,->](2.05,2.55)--(2.5,3);
\draw[thick](1.95,1.55)--(1.5,2)--(1.95,2.45);
\draw[thick](2.5,1)--(2.05,1.45);
\end{tikzpicture}
\qquad
\begin{tikzpicture}[scale=0.95]
\draw[thick](0,0)--(0.9,0.9);
\draw[thick](1.1,1.1)--(1.4,1.4);
\draw[thick,->](1.6,1.6)--(2,2);
\draw[thick,<-](0,1.5)--(2,1.5);
\draw[thick](0,2)--(0.4,1.6);
\draw[thick,->](0.6,1.4)--(2,0);
\draw[thick,<->](2.25,1)--(2.75,1);
\node at(2.5,-1){$\Omega_{3a}$};
\draw[thick,<-](3,0.5)--(5,0.5);
\draw[thick](3,0)--(3.4,0.4);
\draw[thick](3.6,0.6)--(3.9,0.9);
\draw[thick,->](4.1,1.1)--(5,2);
\draw[thick](3,2)--(4.4,0.6);
\draw[thick,->](4.6,0.4)--(5,0);
\end{tikzpicture}
\caption{A generating set of the Reidemeister moves.}
\label{generating}
\end{figure}

In order to show that $F_D(u,v)$ is an invariant for knotoids, we just need to show that it is invariable under the Reidemeister moves.
In \cite{15}, Polyak came up with the minimal generating set of all possible oriented Reidemeister moves for classical knot diagrams. Since the Reidemeister moves for knotoids are the same as the Reidemeister moves for knots, the subset of the oriented Reidemeister moves as shown in Figure~\ref{generating} can also be considered as a generating set of the Reidemeister moves for knotoid diagrams.

\begin{figure}
  \centering
 \subfigure[$\Omega_{1a}$] {\begin{tikzpicture}[scale=1.2]
\draw[dashed](0,0) circle [radius=1];
\fill (0.866,0.5) circle [radius=0.075];
\fill (0.866,-0.5) circle[radius=0.075];
\draw[thick](-0.5,0.866) arc (120:240:1);
\draw[thick,<->](1.5,0)--(2,0);
\draw[dashed](3.5,0) circle [radius=1];
\fill (4.366,0.5) circle [radius=0.075];
\fill (4.366,-0.5) circle[radius=0.075];
\draw[ thick](3,0.866) arc (120:240:1);
\draw[thick,->](2.634,0.5)--(2.634,-0.5);
\node [above]at(2.634,0.5){$+$};
\node [below]at(2.634,-0.5){$-$};
\node[right] at(2.634,0.2){$c'$};
\label{1}
\end{tikzpicture}}
\qquad

 \subfigure[$\Omega_{1b}$]{ \begin{tikzpicture}[scale=1.2]
\draw[dashed](0,0) circle [radius=1];
\fill (0.866,0.5) circle [radius=0.075];
\fill (0.866,-0.5) circle[radius=0.075];
\draw[thick](-0.5,0.866) arc (120:240:1);
\draw[thick,<->](1.5,0)--(2,0);
\draw[dashed](3.5,0) circle [radius=1];
\fill (4.366,0.5) circle [radius=0.075];
\fill (4.366,-0.5) circle[radius=0.075];
\draw[thick](3,0.866) arc (120:240:1);
\draw[thick,<-](2.634,0.5)--(2.634,-0.5);
\node [above]at(2.634,0.5){$-$};
\node [below]at(2.634,-0.5){$$+$$};
\node [right]at(2.634,0.2){$c'$};
\label{2}
\end{tikzpicture}}
\qquad
 \subfigure[$\Omega_{2a}$] {\begin{tikzpicture}[scale=1.2]
\draw[dashed](-3.5,0) circle [radius=1];
\fill (-2.634,0.5) circle [radius=0.075];
\fill (-2.634,-0.5) circle[radius=0.075];
\draw[thick](-3,0.866) arc (60:120:1);
\draw[thick](-4,-0.866) arc (240:300:1);
\draw[thick,<->](-2,0)--(-1.5,0);
\draw[dashed](0,0) circle [radius=1];
\fill (0.866,0.5) circle [radius=0.075];
\fill (0.866,-0.5) circle[radius=0.075];
\draw[thick](0.5,0.866) arc (60:120:1);
\draw[thick](-0.5,-0.866) arc (240:300:1);
\draw[thick,->](0.4,0.9)--(-0.4,-0.9);
\draw[thick,->](-0.4,0.9)--(0.4,-0.9);
\node [above]at(-0.4,0.9){$-$};
\node [above]at(0.4,0.9){+};
\node [above]at(-0.45,-0.9){$-$};
\node [above]at(0.45,-0.9){+};
\node [below]at(0.45,0.7){$c_1$};
\node [below]at(-0.4,0.7){$c_2$};
\label{3}
\end{tikzpicture}}
\qquad
\\
 \subfigure[$\Omega_{3a}$]{ \begin{tikzpicture}[scale=1.2]
\draw[dashed](0,0) circle [radius=1];
\draw[thick](0.8192,0.5736) arc (35:75:1);
\draw[thick](-0.2588,0.9659) arc (105:145:1);
\draw[thick](-0.2588,-0.9659) arc (255:285:1);
\draw[thick,<-](0.34202,0.9397)--(-0.08712,-0.9962);
\draw[thick,<-](-0.34202,0.9397)--(0.08712,-0.9962);
\draw[thick,<-](0.76,0.6428)--(-0.766,0.6428);
\draw[thick,<->](1.5,0)--(2,0);
\draw [dashed](3.5,0) circle [radius=1];
\draw[thick](4.2071,0.7071) arc (45:80:1);
\draw[thick](3.2412,0.9659) arc (105:140:1);
\draw[thick](3.2412,-0.9659) arc (255:290:1);
\draw[thick,<-](3.86,0.9063)--(3.0773,0.9063);
\draw[thick,->](3.4128,-0.996)--(2.9264,0.8192);
\draw[thick,<-](4.0736,0.8192)--(3.6736,-0.9848);
\node [right]at(3.7,1){$-$};
\node [left]at(3.2,1){$+$};
\node [below]at(3.7,-0.9848){$-$};
\node [right]at(4,0.9){$+$};
\node [left]at(3,0.8){$-$};
\node [below]at(3.4128,-0.996){$+$};
\node [right]at(0.6,0.766){$-$};
\node [left]at(-0.6,0.766){$+$};
\node [left]at(-0.3,1){$-$};
\node [below]at(0.1736,-0.9397){$+$};
\node [below]at(-0.1736,-0.9397){$-$};
\node [right]at(0.3,1){$+$};
\node [below]at(0.45,0.7){$c_2$};
\node [below]at(-0.4,0.7){$c_1$};
\node [above]at(0,0.58){$c_3$};
\node [below]at(2.8,0.5){$c'_1$};
\node [below]at(4.2,0.5){$c'_2$};
\node [below]at(3.5,1){$c'_3$};
\fill (0.9659,0.2569) circle [radius=0.075];
\fill (0.9659,-0.2569) circle[radius=0.075];
\fill (4.4659,0.2569) circle [radius=0.075];
\fill (4.4659,-0.2569) circle[radius=0.075];
\label{4}
\end{tikzpicture}}
\qquad
 \subfigure[$\Omega_{3a'}$] {\begin{tikzpicture}[scale=1.2]
\draw[dashed](0,0) circle [radius=1];
\fill (0.9659,0.2569) circle [radius=0.075];
\fill (0.9659,-0.2569) circle[radius=0.075];
\draw[thick](0.766,0.6428) arc (40:80:1);
\draw[thick](-0.2588,0.9659) arc (105:135:1);
\draw[thick](-0.2588,-0.9659) arc (255:290:1);
\draw[thick,<-](0.4226,0.9063)--(-0.4226,0.9063);
\draw[thick,<-](-0.0872,-0.996)--(-0.573576,0.81915);
\draw[thick,<-](0.64278,0.766)--(0.1736,-0.9848);
\draw[thick,<->](1.5,0)--(2,0);
\draw [dashed](3.5,0) circle [radius=1];
\fill (4.4659,0.2569) circle [radius=0.075];
\fill (4.4659,-0.2569) circle[radius=0.075];
\draw[thick](4.266,0.6428) arc (40:75:1);
\draw[thick](3.2412,0.9659) arc (105:140:1);
\draw[thick](3.2412,-0.9659) arc (255:285:1);
\draw[thick,<-](3.84202,0.9397)--(3.3264,-0.9848);
\draw[thick,->](3.15798,0.9397)--(3.6736,-0.9848);
\draw[thick,<-](4.2017,0.7071)--(2.7929,0.7071);
\node [right]at(3.8,1){$-$};
\node [left]at(3.15798,0.9397){$+$};
\node [below]at(3.6736,-0.9848){$-$};
\node [right]at(4.1,0.766){$+$};
\node [left]at(2.9,0.7){$-$};
\node [below]at(3.38,-0.996){$+$};
\node [right]at(0.6,0.766){$-$};
\node [left]at(-0.6,0.766){$+$};
\node [left]at(-0.34202,1){$-$};
\node [below]at(0.1736,-0.9397){$+$};
\node [below]at(-0.1736,-0.9397){$-$};
\node [right]at(0.34202,1){$+$};
\node [below]at(0.3,0.6){$c_3$};
\node [below]at(-0.2,0.6){$c_1$};
\node [above]at(0,0.5){$c_2$};
\node [below]at(3.15,0.5){$c'_1$};
\node [below]at(3.95,0.5){$c'_3$};
\node [above]at(3.5,0.6){$c'_2$};
\end{tikzpicture}}
  \caption{The Gauss diagrams of the Reidemeister moves generating set.}
  \label{gauss diagram generating}
\end{figure}
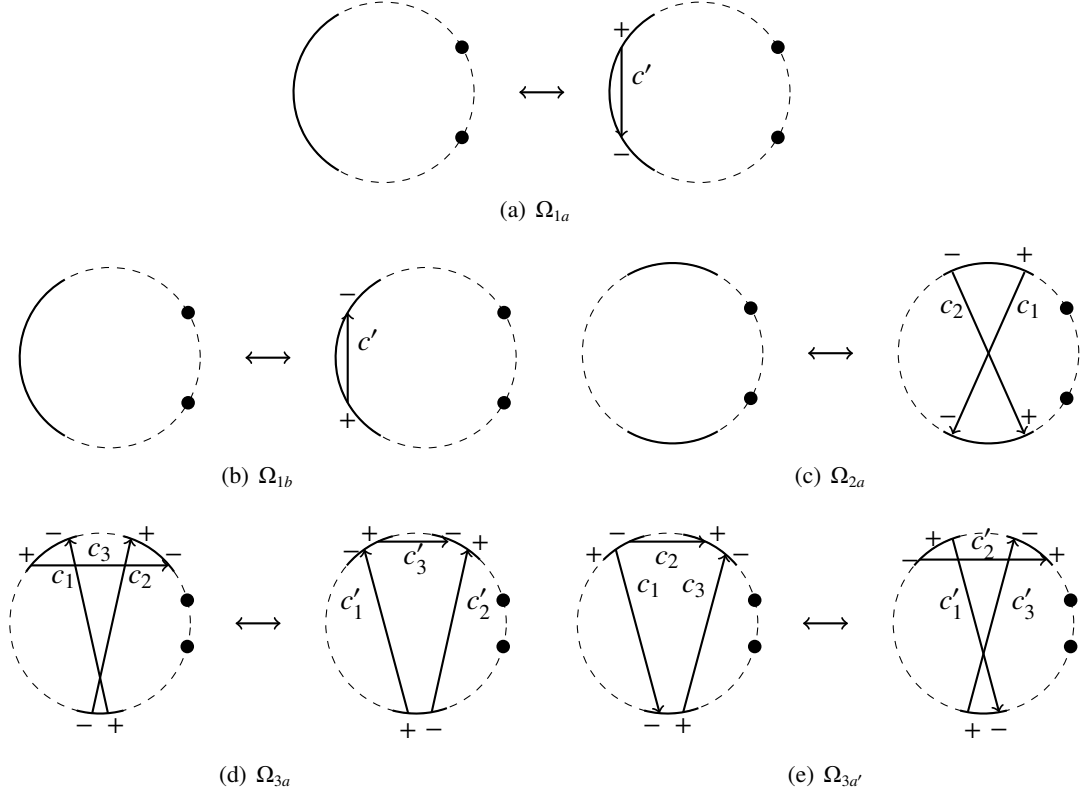

The $F$-polynomial is defined from Gauss diagram, so we mainly consider the Gauss diagrams of the Reidemeister moves {$\Omega_{1a}$}, {$\Omega_{1b}$}, {$\Omega_{2a}$}, {$\Omega_{3a}$} and {$\Omega_{3a'}$} as shown in Figure~\ref{gauss diagram generating}.

\begin{lemma}
Let $D$ be a knotoid diagram. Then $F_D(u,v)$ is an invariant under the moves {$\Omega_{1a}$} and {$\Omega_{1b}$}.
\label{lem1}
\end{lemma}
\begin{proof}
Let $D$ be a knotoid diagram and $D'$ the knotoid diagram obtained from $D$ by once {$\Omega_{1a}$}-move or {$\Omega_{1b}$}-move. Suppose $G(D)$ and $G(D')$ are the corresponding Gauss diagrams of $D$ and $D'$ in Figures~\ref{gauss diagram generating}~\subref{1} and~\ref{gauss diagram generating}~\subref{2}.
Without loss of generality,  we assume that the number of crossings of $D'$ is more than the number of crossings of $D$, and $c'$ is the new crossing in $D'$. Because $c'$ is isolated, $g_{c'}(v)=0$. Then
$$
\begin{aligned}
  F_{D'}(u,v)&=\sum_{c\in C(G(D'))\backslash\left\{c'\right\}}{\rm sign}(c)(u^{g_{c}(v)}-1)+{\rm sign}(c')(u^{g_{c'}(v)}-1)\\&
  =\sum_{c\in C(G(D))}{\rm sign}(c)(u^{g_{c}(v)}-1)=F_{D}(u,v).
  \end{aligned}
$$
This completes the proof.
\end{proof}
\begin{lemma}
Let $D$ be a knotoid diagram. Then $F_D(u,v)$ is an invariant under the move {$\Omega_{2a}$}.
\label{lem2}
\end{lemma}
\begin{proof}
Let $D$ be a knotoid diagram and $D'$ the knotoid diagram obtained from $D$  by once {$\Omega_{2a}$}-move.
Suppose $G(D)$ and $G(D')$ are the corresponding Gauss diagrams of $D$ and $D'$ in Figure~\ref{gauss diagram generating}~\subref{3}. Without loss of generality, we assume that the number of crossings of $D'$ is more than the number of crossings of $D$, $c_1$ and $c_2$ are the new crossings in $D'$. Then ${\rm sign}(c_1)=-{\rm sign}(c_2)$. All chords that pass through $c_1$  are the same as the chords that pass through $c_2$. So we have $l(c_1) = l(c_2)$, $r(c_1) = r(c_2)$ and $g_{c_1}(v)=g_{c_2}(v)$, then
$$
\begin{aligned}
F_{D'}(u,v)&=\sum_{c\in C(G(D'))\backslash\left\{c_1,c_2\right\}}{\rm sign}(c)(u^{g_{c}(v)}-1) + {\rm sign}(c_1)(u^{g_{c_1}(v)}-1)\\&
\ \ \ \ +{\rm sign}(c_2)(u^{g_{c_2}(v)}-1)\\&
=\sum_{c\in C(G(D))}{\rm sign}(c)(u^{g_{c}(v)}-1)=F_{D}(u,v).
\end{aligned}
$$
This completes the proof.
\end{proof}
\begin{lemma}
Let $D$ be a knotoid diagram. Then $F_D(u,v)$ is an invariant under the moves {$\Omega_{3a}$} and {$\Omega_{3a'}$}.
\label{lem3}
\end{lemma}
\begin{proof}
Let $D$ be a knotoid diagram and $D'$ the knotoid diagram obtained from $D$  by once {$\Omega_{3a}$}-move or {$\Omega_{3a'}$}-move.
Suppose $G(D)$ and $G(D')$ are the corresponding Gauss diagrams of $D$ and $D'$ in Figure~\ref{RIII}. Without loss of generality, we assume that  $c_1$, $c_2$ and $c_3$ are crossings in $D$ and $c'_1,~c'_2$ and $c'_3$ are the corresponding crossings in $D'$.
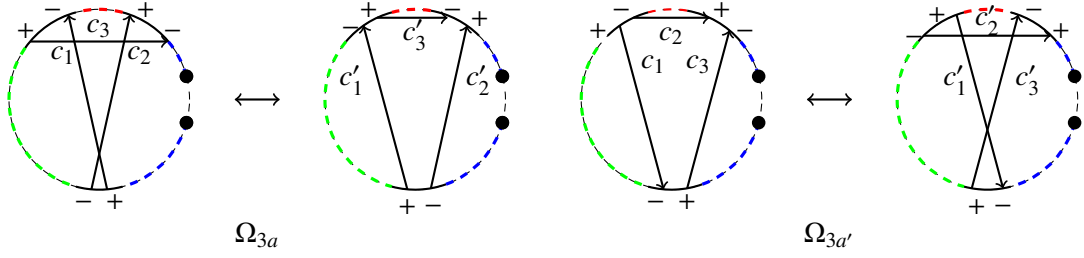
\begin{figure}
  \centering
  \begin{tikzpicture}[scale=1.2]
\draw[dashed](0,0) circle [radius=1];
\draw[green,very thick,dashed](-0.81915,0.573576) arc (145:255:1);
\draw[thick](0.8192,0.5736) arc (35:75:1);
\draw[thick](-0.2588,0.9659) arc (105:145:1);
\draw[thick](-0.2588,-0.9659) arc (255:285:1);
\draw[thick,<-](0.34202,0.9397)--(-0.08712,-0.9962);
\draw[thick,<-](-0.34202,0.9397)--(0.08712,-0.9962);
\draw[thick,<-](0.76,0.6428)--(-0.766,0.6428);
\draw[thick,<->](1.5,0)--(2,0);
\draw [dashed](3.5,0) circle [radius=1];
\draw[green,very thick,dashed](2.7339556,0.6427876) arc (140:255:1);
\draw[thick](4.2071,0.7071) arc (45:80:1);
\draw[thick](3.2412,0.9659) arc (105:140:1);
\draw[thick](3.2412,-0.9659) arc (255:290:1);
\draw[thick,<-](3.86,0.9063)--(3.0773,0.9063);
\draw[thick,->](3.4128,-0.996)--(2.9264,0.8192);
\draw[thick,<-](4.0736,0.8192)--(3.6736,-0.9848);
\node [right]at(3.7,1){$-$};
\node [left]at(3.2,1){$+$};
\node [below]at(3.7,-0.9848){$-$};
\node [right]at(4,0.9){$+$};
\node [left]at(3,0.8){$-$};
\node [below]at(3.4128,-0.996){$+$};
\node [right]at(0.6,0.766){$-$};
\node [left]at(-0.6,0.766){$+$};
\node [left]at(-0.3,1){$-$};
\node [below]at(0.1736,-0.9397){$+$};
\node [below]at(-0.1736,-0.9397){$-$};
\node [right]at(0.3,1){$+$};
\node [below]at(0.45,0.7){$c_2$};
\node [below]at(-0.4,0.7){$c_1$};
\node [above]at(0,0.58){$c_3$};
\node [below]at(2.8,0.5){$c'_1$};
\node [below]at(4.2,0.5){$c'_2$};
\node [below]at(3.5,1){$c'_3$};
\fill (0.9659,0.2569) circle [radius=0.075];
\fill (0.9659,-0.2569) circle[radius=0.075];
\fill (4.4659,0.2569) circle [radius=0.075];
\fill (4.4659,-0.2569) circle[radius=0.075];
\draw[red,very thick,dashed](0.2588,0.9659) arc (75:105:1);
\draw[blue,very thick,dashed](0.342,-0.9397) arc (290:345:1);
\draw[blue,very thick,dashed](0.9659,0.2569) arc (15:40:1);
\draw[red,very thick,dashed](3.6736,0.9848) arc (80:105:1);

\draw[blue,very thick,dashed](3.842,-0.9397) arc (290:345:1);
\draw[blue,very thick,dashed](4.4659,0.2569) arc (15:40:1);
\node at(1.75,-1.5){$\Omega_{3a}$};
\fill (0.9659,0.2569) circle [radius=0.075];
\fill (0.9659,-0.2569) circle[radius=0.075];
\fill (4.4659,0.2569) circle [radius=0.075];
\fill (4.4659,-0.2569) circle[radius=0.075];
\end{tikzpicture}
\qquad
\begin{tikzpicture}[scale=1.2]
\draw[dashed](0,0) circle [radius=1];

\draw[thick](0.766,0.6428) arc (40:80:1);
\draw[thick](-0.2588,0.9659) arc (105:135:1);
\draw[thick](-0.2588,-0.9659) arc (255:290:1);
\draw[thick,<-](0.4226,0.9063)--(-0.4226,0.9063);
\draw[thick,<-](-0.0872,-0.996)--(-0.573576,0.81915);
\draw[thick,<-](0.64278,0.766)--(0.1736,-0.9848);
\draw[thick,<->](1.5,0)--(2,0);
\draw [dashed](3.5,0) circle [radius=1];

\draw[thick](4.266,0.6428) arc (40:75:1);
\draw[thick](3.2412,0.9659) arc (105:140:1);
\draw[thick](3.2412,-0.9659) arc (255:285:1);
\draw[thick,<-](3.84202,0.9397)--(3.3264,-0.9848);
\draw[thick,->](3.15798,0.9397)--(3.6736,-0.9848);
\draw[thick,<-](4.2017,0.7071)--(2.7929,0.7071);
\node at(1.75,-1.5){$\Omega_{3a'}$};
\node [right]at(3.8,1){$-$};
\node [left]at(3.15798,0.9397){$+$};
\node [below]at(3.6736,-0.9848){$-$};
\node [right]at(4.1,0.766){$+$};
\node [left]at(2.9,0.7){$-$};
\node [below]at(3.38,-0.996){$+$};
\node [right]at(0.6,0.766){$-$};
\node [left]at(-0.6,0.766){$+$};
\node [left]at(-0.34202,1){$-$};
\node [below]at(0.1736,-0.9397){$+$};
\node [below]at(-0.1736,-0.9397){$-$};
\node [right]at(0.34202,1){$+$};
\node [below]at(0.3,0.6){$c_3$};
\node [below]at(-0.2,0.6){$c_1$};
\node [above]at(0,0.5){$c_2$};
\node [below]at(3.15,0.5){$c'_1$};
\node [below]at(3.95,0.5){$c'_3$};
\node [above]at(3.5,0.6){$c'_2$};
\draw[red,thick,dashed](0.1736,0.9848) arc (80:105:1);
\draw[green,very thick,dashed](-0.8192,0.5735) arc (145:255:1);
\draw[blue,very thick,dashed](0.342,-0.9397) arc (290:345:1);
\draw[blue,very thick,dashed](0.9659,0.2569) arc (15:40:1);
\draw[red,very thick,dashed](3.6736,0.9848) arc (80:105:1);
\draw[green,very thick,dashed](2.7339,0.6428) arc (140:255:1);
\draw[blue,very thick,dashed](3.842,-0.9397) arc (290:345:1);
\draw[blue,very thick,dashed](4.4659,0.2569) arc (15:40:1);
\fill (4.4659,0.2569) circle [radius=0.075];
\fill (4.4659,-0.2569) circle[radius=0.075];
\fill (0.9659,0.2569) circle [radius=0.075];
\fill (0.9659,-0.2569) circle[radius=0.075];
\end{tikzpicture}
  \caption{{$\Omega_{3a}$} and {{$\Omega_{3a'}$}}.}
  \label{RIII}
\end{figure}

Let $D'$ be the knotoid diagram obtained from $D$  by once {$\Omega_{3a}$}-move. Obviously, ${\rm sign}(c_1)={\rm sign}(c'_1)=+1$, ${\rm sign}(c_2)={\rm sign}(c'_2)=-1$ and ${\rm sign}(c_3)={\rm sign}(c'_3)=+1$. Suppose the sum of the signs of the endpoints on the green arc is $x$, the sum of the signs of the endpoints on the red arc is $y$ and the sum of the signs of the endpoints on the blue arcs is $z$. Then $x+y+z=0$, $i(c_1)=y+z$, $i(c_2)=z$ and $i(c_3)=x+z$. Therefore,
\begin{equation}\label{ieuqation}
i(c_1)+i(c_3)-i(c_2)=x+y+z=0.
\end{equation}
According to the definition of $F$-polynomial, we have
\begin{equation}\label{8}
\begin{aligned}
 F_{D}(u,v)=&\sum_{c\in C(G(D))\backslash\left\{c_1,c_2,c_3\right\} }{\rm sign}(c)(u^{g_c(v)}-1)+{\rm sign}(c_1)(u^{g_{c_1}(v)}-1)\\&
 +{\rm sign}(c_2)(u^{g_{c_2}(v)}-1)+{\rm sign}(c_3)(u^{g_{c_3}(v)}-1),
 \end{aligned}
 \end{equation}
 and
 \begin{equation}\label{12}
\begin{aligned}
 F_{D'}(u,v)=&\sum_{c'\in C(G(D'))\backslash\left\{c'_1,c'_2,c'_3\right\} }{\rm sign}(c')(u^{g_{c'}(v)}-1)+{\rm sign}(c'_1)(u^{g_{c'_1}(v)}-1)\\&
 +{\rm sign}(c'_2)(u^{g_{c'_2}(v)}-1)+{\rm sign}(c'_3)(u^{g_{c'_3}(v)}-1),
\end{aligned}
 \end{equation}
 where
$$
 \begin{aligned}
&g_{c_1}(v)=g_{c'_1}(v)+ v^{\phi_{c_1}(-i(c_2))}-v^{\phi_{c_1}(-i(c_3))},\\&
g_{c_2}(v)=g_{c'_2}(v)+ v^{\phi_{c_2}(i(c_1))}-v^{\phi_{c_2}(-i(c_3))},\\&
g_{c_3}(v)=g_{c'_3}(v)+ v^{\phi_{c_3}(i(c_1))}-v^{\phi_{c_3}(i(c_2))}.
\end{aligned}
$$
According to Equation~(\ref{ieuqation}), we have
$$
\begin{aligned}
&v^{\phi_{c_1}(-i(c_2))}-v^{\phi_{c_1}(-i(c_3))}=v^{\phi_{c_1}(-i(c_1)-i(c_3))}-v^{\phi_{c_1}(-i(c_3))}=0,\\&
v^{\phi_{c_2}(i(c_1))}-v^{\phi_{c_2}(-i(c_3))}=v^{\phi_{c_2}(i(c_2)-i(c_3))}-v^{\phi_{c_2}(-i(c_3))}=0,\\&
v^{\phi_{c_3}(i(c_1))}-v^{\phi_{c_3}(i(c_2))}=v^{\phi_{c_3}(i(c_2)-i(c_3))}-v^{\phi_{c_3}(i(c_2))}=0.
\end{aligned}
$$
Therefore, $g_{c_1}(v)=g_{c'_1}(v)$, $g_{c_2}(v)=g_{c'_2}(v)$ and $g_{c_3}(v)=g_{c'_3}(v)$. According to Equations~(\ref{8}) and (\ref{12}), $ F_{D'}(u,v)= F_{D}(u,v)$.

Let $D'$ be the knotoid diagram obtained from $D$  by once {$\Omega_{3a'}$}-move. Following the lines as above, we have the same conclusion. This completes the proof.
\end{proof}

\begin{proof}[Proof of Theorem~\ref{invariance theorem}]
  According to Lemmas~\ref{lem1}, \ref{lem2} and \ref{lem3}, $F$-polynomial is an invariant for knotoids.

  This completes the proof of Theorem~\ref{invariance theorem}.
\end{proof}

\subsection{The properties of $F$-polynomial}
Next, we mainly consider the properties of $F$-polynomial.

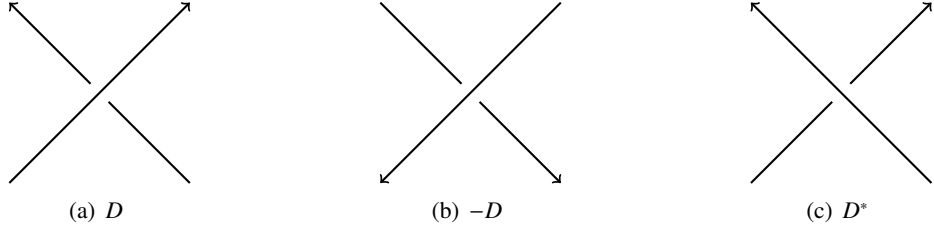
\begin{figure}
\centering
\subfigure[$D$]{
\begin{tikzpicture}[scale=1.2]
\draw[thick,->](0,0)--(2,2);
\draw[thick](2,0)--(1.1,0.9);
\draw[thick,->](0.9,1.1)--(0,2);
\label{21}
\end{tikzpicture}}
\qquad
\qquad
\qquad
\subfigure[$-D$]{
\begin{tikzpicture}[scale=1.2]
\draw[thick,<-](0,0)--(2,2);
\draw[thick,<-](2,0)--(1.1,0.9);
\draw[thick](0.9,1.1)--(0,2);
\label{22}
\end{tikzpicture}}
\qquad
\qquad
\qquad
\subfigure[$D^*$]{
\begin{tikzpicture}[scale=1.2]
\draw[thick,->](2,0)--(0,2);
\draw[thick](0,0)--(0.9,0.9);
\draw[thick,->](1.1,1.1)--(2,2);
\label{23}
\end{tikzpicture}}
\caption{The inverse of $D$ and the mirror image of $D$.}
\label{6}
\end{figure}

For a knotoid diagram $D$, let $-D$ be the knotoid diagram  obtained from $D$ by changing the orientation of $D$. $-D$ is called the $inverse$ of $D$ as shown in Figure~\ref{6}~\subref{22}. Let $D^*$ be the knotoid diagram obtained from $D$ by changing the over/under information of $D$. $D^*$ is called the $mirror ~image$ of $D$ as shown in Figure~\ref{6}~\subref{23}. If $D$ and $D^*$ are equivalent, then $D$ is called $amphicheiral$ and otherwise $cheiral$.
\begin{proposition}\label{pro1}
Let $D$ be a knotoid diagram and  $-D$ the inverse of $D$, then $F_{-D}(u,v)=F_{D}(u,v)$.
\end{proposition}
\begin{proof}
  Let $c$ be the crossing in $D$ and $c'$ the corresponding  crossing  in $-D$. Suppose $G(D)$ and $G(-D)$ are the  Gauss diagrams of $D$ and $-D$ as shown in Figure~\ref{ D and $-$D}. Then ${\rm sign}(c')={\rm sign}(c)$, $i(c')=i(c)$, $r(c')=r(c)$ and $l(c')=l(c)$. So $\phi_{c'}=\phi_c$  and  $g_{c'}(v)=g_{c}(v)$.
  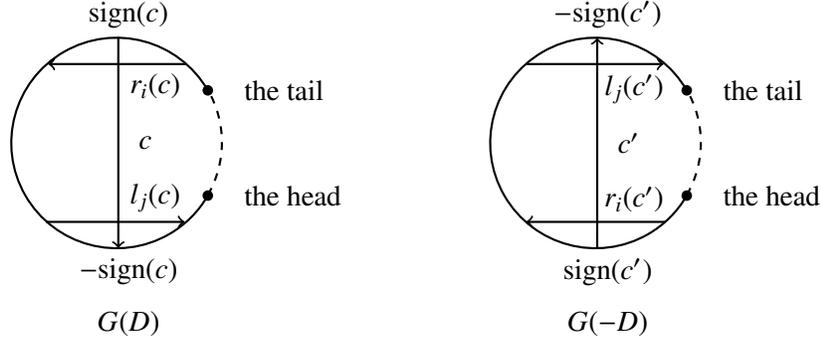
\begin{figure}
    \centering
\begin{tikzpicture}[scale=0.7]
    \draw[thick](0,0) arc (30:330:2);
    \draw[thick,dashed] (0,-1.99) arc (330:390:2);
    \fill (0,0) circle (.1);
    \node at (0.5,0) [right]{the tail};
    \fill (0,-2) circle (.1);
    \node at (0.5,-2) [right]{the head};
    \draw [thick,->](-1.7,1) -- (-1.7,-3);
    \node [above] at(-1.5,1){sign($c$)};
    \node [below]at(-1.5,-3){$-$sign($c$)};
     \node [below]at(-1.5,-4){$G(D)$};
    \node[right] at (-1.5,-1){$c$};
    \draw [thick,->](-0.4,0.5) -- (-3.03,0.5);
    \draw [thick,<-](-0.44,-2.5) -- (-3.075,-2.5);
     \node [below]at(-1,0.5){$r_i(c)$};
     \node [above] at(-1,-2.5){$l_j(c)$};
\end{tikzpicture}
\qquad
\qquad
\begin{tikzpicture}[scale=0.7]
    \draw[thick](0,0) arc (30:330:2);
    \draw[thick,dashed] (0,-1.99) arc (330:390:2);
    \fill (0,0) circle (.1);
    \node at (0.5,0) [right]{the tail};
    \fill (0,-2) circle (.1);
    \node at (0.5,-2) [right]{the head};
    \draw [thick,<-](-1.7,1) -- (-1.7,-3);
    \node [above] at(-1.5,1){$-$sign($c'$)};
    \node [below]at(-1.5,-3){sign($c'$)};
     \node [below]at(-1.5,-4){$G(-D)$};
    \node[right] at (-1.5,-1){$c'$};
     \draw [thick,<-](-0.44,0.5) -- (-3.05,0.5);
    \draw [thick,->](-0.4,-2.5) -- (-3.05,-2.5);
    \node [below]at(-1,0.5){$l_j(c')$};
     \node [above] at(-1,-2.5){$r_i(c')$};
\end{tikzpicture}
    \caption{Gauss diagrams of $D$ and $-D$.}\label{ D and $-$D}
  \end{figure}
 Then
$ \begin{aligned}
 F_{-D}(u,v)
 =F_D(u,v).
 \end{aligned}$
\end{proof}

\begin{proposition}\label{pro2}
Let $D$ be a knotoid diagram and $D^*$  the mirror image of $D$, then $F_{D^*}(u,v)=-F_{D}(u,v^{-1})$.
\end{proposition}
\begin{proof}
  Let $c$ be the crossing in $D$ and $c'$ the corresponding crossing in $D^*$. Suppose $G(D)$ and $G(D^*)$ are the  Gauss diagrams of $D$ and $D^*$ as shown in Figure~\ref{9}. Then ${\rm sign}(c')=-{\rm sign}(c)$, $i(c')=i(c)$, $r(c')=l(c)$ and $l(c')=r(c)$. So $\phi_{c'}=\phi_c$. We consider the relationship between $g_{c'}(v)$ and $g_{c}(v)$.
  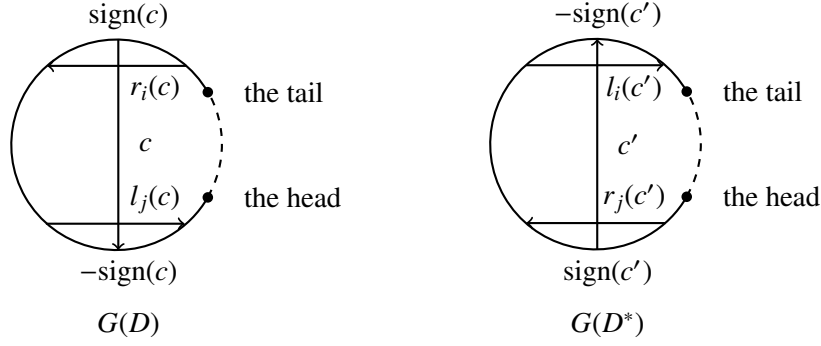
\begin{figure}
    \centering
\begin{tikzpicture}[scale=0.7]
    \draw[thick](0,0) arc (30:330:2);
    \draw[thick,dashed] (0,-1.99) arc (330:390:2);
    \fill (0,0) circle (.1);
    \node at (0.5,0) [right]{the tail};
    \fill (0,-2) circle (.1);
    \node at (0.5,-2) [right]{the head};
    \draw [thick,->](-1.7,1) -- (-1.7,-3);
    \node [above] at(-1.5,1){sign($c$)};
    \node [below]at(-1.5,-3){$-$sign($c$)};
     \node [below]at(-1.5,-4){$G(D)$};
    \node[right] at (-1.5,-1){$c$};
    \draw [thick,->](-0.4,0.5) -- (-3.05,0.5);
    \draw [thick,<-](-0.44,-2.5) -- (-3.075,-2.5);
     \node [below]at(-1,0.5){$r_i(c)$};
     \node [above] at(-1,-2.5){$l_j(c)$};
\end{tikzpicture}
\qquad
\qquad
\begin{tikzpicture}[scale=0.7]
    \draw[thick](0,0) arc (30:330:2);
    \draw[thick,dashed] (0,-1.99) arc (330:390:2);
    \fill (0,0) circle (.1);
    \node at (0.5,0) [right]{the tail};
    \fill (0,-2) circle (.1);
    \node at (0.5,-2) [right]{the head};
    \draw [thick,<-](-1.7,1) -- (-1.7,-3);
    \node [above] at(-1.5,1){$-$sign($c'$)};
    \node [below]at(-1.5,-3){sign($c'$)};
     \node [below]at(-1.5,-4){$G(D^*)$};
    \node[right] at (-1.5,-1){$c'$};
     \draw [thick,<-](-0.44,0.5) -- (-3.05,0.5);
    \draw [thick,->](-0.4,-2.5) -- (-3.05,-2.5);
    \node [below]at(-1,0.5){$l_i(c')$};
     \node [above] at(-1,-2.5){$r_j(c')$};
\end{tikzpicture}
    \caption{Gauss diagrams of $D$ and $D^*$.}
    \label{9}
  \end{figure}

$$
  g_c(v)=\sum^n_{i=1}{\rm sign}(r_i(c))v^{\phi_c(i(r_i(c)))}-\sum^m_{j=1}{\rm sign}(l_j(c))v^{\phi_c(-i(l_j(c)))},
$$
and
$$
    \begin{aligned}
  g_{c'}(v)&=\sum^m_{i=1}{\rm sign}(r_i(c'))v^{\phi_{c'}(i(r_i(c')))}-\sum^n_{j=1}{\rm sign}(l_j(c'))v^{\phi_{c'}(-i(l_j(c')))}
 \\&= -\sum^m_{i=1}{\rm sign}(l_i(c))v^{\phi_{c}(i(l_i(c)))}+\sum^n_{j=1}{\rm sign}(r_j(c))v^{\phi_{c}(-i(r_j(c)))}
\\&=\sum^n_{i=1}{\rm sign}(r_i(c))v^{-\phi_c(i(r_i(c)))}-\sum^m_{j=1}{\rm sign}(l_j(c))v^{-\phi_c(-i(l_j(c)))}
\\&= g_c(v^{-1}).
 \end{aligned}
$$
 Then
$$
 \begin{aligned}
 F_{D^*}(u,v)&=\sum_{c'\in C(G(D^*))}{\rm sign}(c')(u^{g_{c'}(v)}-1)\\&
 =-\sum_{c\in C(G(D))}{\rm sign}(c)(u^{g_c(v^{-1})}-1)\\&
 =-F_D(u,v^{-1}).
 \end{aligned}
$$
 This completes the proof.
\end{proof}

Therefore, for a given knotoid diagram, $F_D(u,v)$ can distinguish whether  it$'$s amphicheiral or not.

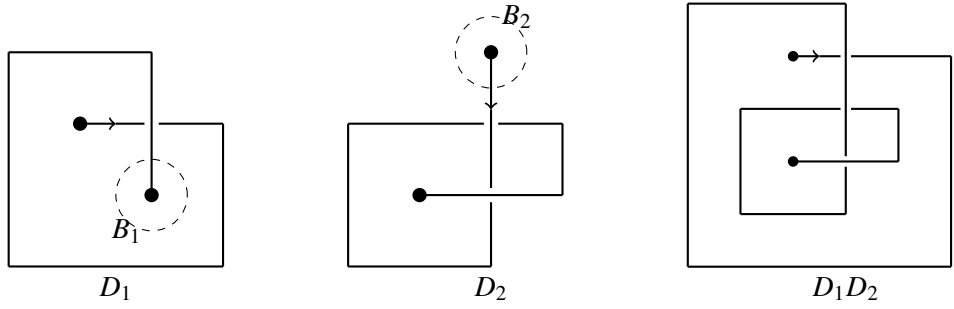
\begin{figure}
\begin{tikzpicture}[scale=0.95]
\draw[thick,->](-1,1)--(-0.5,1);
\draw[thick](-0.5,1)--(-0.1,1);
\draw[thick](0.1,1)--(1,1);
\draw[thick](1,1)--(1,-1);
\draw[thick](1,-1)--(-2,-1);
\draw[thick](-2,-1)--(-2,2);
\draw[thick](-2,2)--(0,2);
\draw[thick](0,2)--(0,0);
\fill (0,0) circle (.1);
\fill (-1,1) circle (.1);
\draw[dashed](0,0) circle [radius=0.5];
\node [below]at(-0.5,-1){$D_1$};
\node [left]at(0,-0.5){$B_1$};
\end{tikzpicture}
\qquad
\qquad
\begin{tikzpicture}[scale=0.95]
\draw[thick,->](0,2)--(0,1.2);
\draw[thick](0,1.2)--(0,0.1);
\draw[thick](0,-0.1)--(0,-1);
\draw[thick](0,-1)--(-2,-1);
\draw[thick](-2,-1)--(-2,1);
\draw[thick](-2,1)--(-0.1,1);
\draw[thick](0.1,1)--(1,1);
\draw[thick](1,1)--(1,0);
\draw[thick](1,0)--(-1,0);
\fill (0,2) circle (.1);
\fill (-1,0) circle (.1);
\draw[dashed](0,2) circle [radius=0.5];
\node [below]at(0,-1){$D_2$};
\node [right]at(0,2.5){$B_2$};
\end{tikzpicture}
\qquad
\qquad
\begin{tikzpicture}[scale=0.7]
\draw[thick,->](-1,1)--(-0.5,1);
\draw[thick](-0.5,1)--(-0.1,1);
\draw[thick](0.1,1)--(2,1);
\draw[thick](2,1)--(2,-3);
\draw[thick](2,-3)--(-3,-3);
\draw[thick](-3,-3)--(-3,2);
\draw[thick](-3,2)--(0,2);
\draw[thick](0,2)--(0,-0.9);
\draw[thick](0,-1.1)--(0,-2);
\draw[thick](0,-2)--(-2,-2);
\draw[thick](-2,-2)--(-2,0);
\draw[thick](-2,0)--(-0.1,0);
\draw[thick](0.1,0)--(1,0);
\draw[thick](1,0)--(1,-1);
\draw[thick](1,-1)--(-1,-1);
\fill (-1,1) circle (.1);
\fill (-1,-1) circle (.1);

\node [below]at(0,-3){$D_1D_2$};

\end{tikzpicture}
\caption{Product of  $D_1$ and $D_2$.}
\label{product}
\end{figure}

Given two knotoids $d_1$ and $d_2$. Take $D_1$ and $D_2$ are the knotoid diagrams to represent them, respectively. For the head of $D_1$ and the tail of $D_2$, pick two disks $B_1$ and $B_2$, such that $D_i$ meets $B_i$ precisely along a radius of $B_i$ and all crossings outside $B_i$ for $i=1,2$. Then glue $S^2-$Int$(B_1)$ and $S^2-$Int$(B_2)$ along the boundary homeomorphism of the disks, carrying the only point of $D_1\cap \partial B_1 $ to the only point of $D_2\cap \partial B_2$. The part of $D_1$ lying in $S^2-$Int$(B_1)$ and the part of $D_2$ lying in $S^2-$Int$(B_2)$ form a new knotoid diagram denoted by  $D_1D_2$, called the $product $ of $D_1$ and $D_2$ as shown in Figure~\ref{product}. The $product$ $d_1d_2$ is the equivalence class of $D_1D_2$. It is well-defined since any knotoid $d$ in  $S^2$ can be represented by a knotoid diagram $D$ that the tail of $D$ is in the exterior region\cite{16}.

\begin{proposition}
Let $D$ and $D'$ be knotoid diagrams, $DD'$  the product of $D$ and $D'$. Then $F_{DD'}(u,v)=F_D(u,v)+F_{D'}(u,v)$.
\end{proposition}
\begin{proof}
Suppose that the tail of $D'$ is in the exterior region. According to the definition of $DD'$, $DD'$ does not produce any new crossing and the original crossings remain unchanged. Therefore, for the Gauss diagram of $DD'$, the part of the chords contributed by $D$ does not intersect with the chords contributed by $D'$ and the original  intersection relationship between the chords remains unchanged. Obviously,  $F_{DD'}(u,v)=F_D(u,v)+F_{D'}(u,v)$.

If the tail of $D'$ is in the interior region, we can find a knotoid diagram $D''$ that is equivalent to $D'$. And the tail of $D''$ is in the exterior region. So, $F_{DD''}(u,v)=F_D(u,v)+F_{D''}(u,v)$. For the part $D''$ of $DD''$, we can deform it into $D'$ by isotopies of $S^2$ and classical Reidemeister moves taking place in a local disk free of the endpoints which do not affect the part $D$ of $DD''$. Therefore, we have $F_{DD'}(u,v)=F_{DD''}(u,v)=F_D(u,v)+F_{D''}(u,v)$. Because $D''$ that is equivalent to $D'$,  $F_{D'}(u,v)=F_{D''}(u,v)$. Then $F_{DD'}(u,v)=F_D(u,v)+F_{D'}(u,v)$.

In conclusion, we have $F_{DD'}(u,v)=F_D(u,v)+F_{D'}(u,v)$.
\end{proof}
\begin{proposition}\label{pro3}
For each knot-type knotoid diagram D, $F_D(u,v)$ = $0$.
\end{proposition}
\begin{proof}
Let $D$ be a knot-type knotoid diagram and $G(D)$ the  Gauss diagram of $D$. For each chord $c$  in $G(D)$, we have $i(c)=0$. So we can get $g_c(v)=0$, and  then $F_D(u,v)=0$.
\end{proof}

Here, we give an example. The inverse of Proposition 3.10 is not true.

\begin{example}
{\rm As shown in Figure~\ref{anti}, $D$ is a knotoid diagram of 5.1.26 of
Bartholomew$'$s table\cite{1}. It is a proper knotoid diagram, but $F_D(u,v)=0$.}
\end{example}
\begin{figure}
  \centering
 \begin{tikzpicture}[scale=0.9]
\draw[thick,->](0,0)--(0,-0.5);
\draw[thick](0,-0.5)--(0,-2);
\draw[thick](0,-2)--(-0.9,-2);
\draw[thick](-1.1,-2)--(-2,-2);
\draw[thick](-2,-2)--(-2,1);
\draw[thick](-2,1)--(-4,1);
\draw[thick](-4,1)--(-4,-1);
\draw[thick](-4,-1)--(-2.1,-1);
\draw[thick](-1.9,-1)--(-0.1,-1);
\draw[thick](0.1,-1)--(1,-1);
\draw[thick](1,-1)--(1,1);
\draw[thick](1,1)--(-1,1);
\draw[thick](-1,1)--(-1,-0.9);
\draw[thick](-1,-1.1)--(-1,-3);
\draw[thick](-1,-3)--(-3,-3);
\draw[thick](-3,-3)--(-3,-1.1);
\draw[thick](-3,-0.9)--(-3,0);
\fill (0,0) circle (.1);
\fill (-3,0) circle (.1);
\node [below] at (0.2,-1){$a$};
\node [below] at (-0.8,-2){$b$};
\node [above] at (-1.8,-1){$c$};
\node [above] at (-0.8,-1){$e$};
\node [below] at (-2.8,-1){$d$};
\end{tikzpicture}
  \caption{A proper knotoid diagram $D$.}
  \label{anti}
\end{figure}
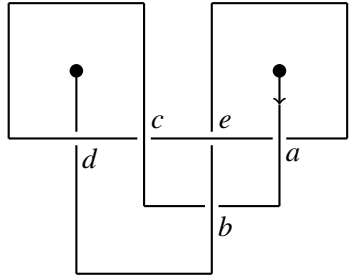
\begin{proposition}
Let $D$ be a knotoid diagram, then $F_D(t,1)=F_D(t)$.
\end{proposition}
\begin{proof}
Let $G(D)$ be the Gauss diagram of $D$. For each chord $c$ in $G(D)$, we consider $g_c(1)$. According to the definitions of $l(c)$ and $r(c)$, we have
$$
g_c(1)=\sum_{i=1}^{n}{\rm sign}(r_i(c))-\sum_{j=1}^{m}{\rm sign}(l_j(c))=i(c).
$$
Let $u=t$ and $v=1$, then we have $$F_D(u,v)=F_D(t,1)=\sum_{c\in C(G(D))}{\rm sign}(c)(t^{g_c(1)}-1)=\sum_{c\in C(G(D))}{\rm sign}(c)(t^{i(c)}-1)=F_D(t).$$
\end{proof}

In this view, $F$-polynomial can be considered as a generalization of the index polynomial.
\subsection{$n$th $F$-polynomial for knotoids}
We can also get a class of polynomials from the $F$-polynomial.

Next, we define the $n$th $F$-polynomial for knotoid diagram for each nonnegative integer $n$.
Let $D$ be a  knotoid diagram, $G(D)$  the  Gauss diagram of $D$ and $C(G(D))$ the set of chords of $G(D)$. For each nonnegative integer  $n$,
we take
$$C_{n}(G(D)) = \left\{
c\in C(G(D))|i(c)=kn~{\rm for~some~integer}~k
\right\}.$$
Then take the $d_n$-function  as the sum of signs
of the endpoints on the lead side of $c$ whose chords in $C_{n}(G(D))\backslash
\left\{
 c
 \right\}$(\cite{11}).

For each chord $c\in C_n(G(D))$, we define two subsets of  $C_{n}(G(D))$. One is the chords that pass through $c$ and point to the lead side of  $c$ denoted by $l^n(c)= \left\{l^n_1(c), l^n_2(c),{\cdots},l^n_q(c) \right\}$ and the other is the chords that pass through $c$ and point from the lead side of  $c$ to the other side denoted by $r^n(c)=\left\{r^n_1(c), r^n_2(c),{\cdots},r^n_p(c) \right\}$.
 For each chord $c\in C_{n}(G(D))$, take a map $\phi^n_c$ from $\mathbb{Z}$ to $\mathbb{Z}_{|d_n(c))|}$,  which maps the $d_n$-functions of the chords that pass through  $c$ to $\mathbb{Z}_{|d_n(c)|}$.
\begin{definition}
{\rm Let $D$ be a  knotoid diagram and $G(D)$ the Gauss diagram of $D$, for each nonnegative integer $n$ and each chord $c$ in $G_n(D)$, the $nth~ index~ function~ of~ c$ is defined by}
\begin{equation}
g^n_c(v)=\sum^p_{i=1}{\rm sign}(r^n_i(c))v^{\phi^n_c(d_n(r^n_i(c)))}-\sum^q_{j=1}{\rm sign}(l^n_j(c))v^{\phi^n_c(-d_n(l^n_j(c)))}.
\label{gcn}
\end{equation}
\end{definition}

Then we give a class of polynomials as follows.
\begin{definition}
{\rm Let $D$ be a knotoid diagram, $G(D)$ the  Gauss diagram of $D$ and  $C_n(G(D))$  the subset of chords of $G(D)$. For each nonnegative interger $n$, the $nth~F$-$polynomial$  $F^n_D(u,v)$ is defined by}
\begin{equation}
F^n_D(u,v)=\sum_{c\in C_n(G(D))}{\rm sign}(c)(u^{g^n_c(v)}-1).
\end{equation}
\end{definition}

\begin{theorem}
For each nonnegative integer $n$,  $F^n_D(u,v)$ is an invariant for  knotoids.
\end{theorem}
\begin{proof}
The proof is similar to the proof of Theorem~\ref{invariance theorem}, we only need to consider the invariance of the Reidemeister moves {$\Omega_{1a}$}, {$\Omega_{1b}$}, {$\Omega_{2a}$}, {$\Omega_{3a}$} and {$\Omega_{3a'}$}.
\end{proof}
\begin{remark}
For $n=1$, $F^1_D(u,v)$ is the $F$-polynomial for  knotoids.
\end{remark}
\begin{remark}
The definition of the $n$th $F$-polynomial for knotoids can be extended to virtual knotoids.
\end{remark}

  $F^n_D(u,v)$ has similar properties as $F_D(u,v)$.
 \begin{proposition}
Let $D$ and $D'$ be knotoid diagrams, $-D$  the inverse of $D$, $D^*$ the mirror image of $D$ and $DD'$ the product of $D$ and $D'$. Then

$\begin{aligned}
&(1)\ F^n_{-D}(u,v)=F^n_{D}(u,v),\\&
(2)\  F^n_{D^*}(u,v)=-F^n_{D}(u,v^{-1}),\\&
 (3)\  F^n_{DD'}(u,v)=F^n_{D}(u,v)+F^n_{D'}(u,v).
 \end{aligned}$
\end{proposition}

\section{A family of knotoid diagrams}\label{A family of knotoid diagrams}
In this section, we will  provide a family of knotoid diagrams which can be distinguished from each other by $F$-polynomial but cannnot be distinguished by the index polynomial and the $n$th  polynomial.

\begin{example}
{\rm As shown in Figure~\ref{D1}, $D_1$ and $D_2$ are the knotoid diagrams of 5.1.3 and 5.1.4 of Bartholo-\\
mew$'$s table\cite{1}, respectively. We calculate $F_{D_i}(t)$, $Z^n_{D_i}(t)$ and $F_{D_i}(u,v)$ ($i=1,2$).}
\end{example}
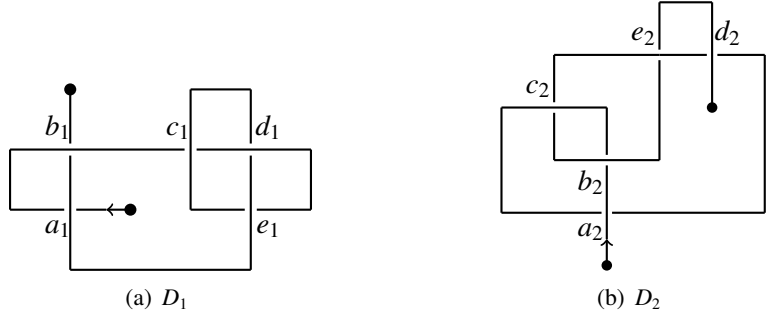
\begin{figure}
  \centering
  \subfigure[$D_1$]{
\begin{tikzpicture}[scale=0.8]
\draw[thick,->](0,0)--(-0.4,0);
\draw[thick](-0.4,0)--(-0.9,0);
\draw[thick](-1.1,0)--(-2,0);
\draw[thick](-2,0)--(-2,1);
\draw[thick](-2,1)--(0.9,1);
\draw[thick](1.1,1)--(3,1);
\draw[thick](3,1)--(3,0);
\draw[thick](3,0)--(2.1,0);
\draw[thick](1.9,0)--(1,0);
\draw[thick](1,0)--(1,2);
\draw[thick](1,2)--(2,2);
\draw[thick](2,2)--(2,1.1);
\draw[thick](2,0.9)--(2,-1);
\draw[thick](2,-1)--(-1,-1);
\draw[thick](-1,-1)--(-1,0.9);
\draw[thick](-1,1.1)--(-1,2);
\fill (0,0) circle (.1);
\fill (-1,2) circle (.1);
\node [below] at (-1.2,0){$a_1$};
\node [above] at (-1.2,1){$b_1$};
\node [above] at (0.8,1){$c_1$};
\node [above] at (2.3,1){$d_1$};
\node [below] at (2.3,0){$e_1$};
\end{tikzpicture}}
\qquad
\qquad
\qquad
\subfigure[$D_2$]{
\begin{tikzpicture}[scale=0.7]
\draw[thick,->](0,0)--(0,0.5);
\draw[thick](0,0.5)--(0,1.9);
\draw[thick](0,2.1)--(0,3);
\draw[thick](0,3)--(-2,3);
\draw[thick](-2,3)--(-2,1);
\draw[thick](-2,1)--(-0.1,1);
\draw[thick](0.1,1)--(3,1);
\draw[thick](3,1)--(3,4);
\draw[thick](3,4)--(2.1,4);
\draw[thick](1.9,4)--(1,4);
\draw[thick](1,4)--(-1,4);
\draw[thick](-1,4)--(-1,3.1);
\draw[thick](-1,2.9)--(-1,2);
\draw[thick](-1,2)--(1,2);
\draw[thick](1,2)--(1,3.9);
\draw[thick](1,4.1)--(1,5);
\draw[thick](1,5)--(2,5);
\draw[thick](2,5)--(2,3);
\fill (0,0) circle (.1);
\fill (2,3) circle (.1);
\node [below] at (-0.3,1){$a_2$};
\node [below] at (-0.3,2){$b_2$};
\node [above] at (-1.3,3){$c_2$};
\node [above] at (0.7,4){$e_2$};
\node [above] at (2.3,4){$d_2$};
\end{tikzpicture}}

  \caption{Two knotoid diagrams.}
  \label{D1}
\end{figure}

\begin{figure}
  \centering
  \subfigure[$G(D_1)$]{
  \begin{tikzpicture}[scale=0.8]
\draw[dashed](0,0) circle [radius=2];
\draw[thick](1.9318,0.5176) arc (15:345:2);
\draw[thick,<-](1.4142,1.4142)--(0.5176,-1.9318);
\draw[thick,->](0.5176,1.9318)--(1.4142,-1.4142);
\draw[thick,<-](-0.5176,1.9318)--(-1.9318,-0.5176);
\draw[thick,->](-1.4142,1.4142)--(-1.4142,-1.4142);
\draw[thick,<-](-1.9318,0.5176)--(-0.5176,-1.9318);
\fill (1.9318,0.5176) circle [radius=0.075];
\fill (1.9318,-0.5176) circle [radius=0.075];
\node [right]at(1.4142,1.4142){$-$};
\node [below]at(0.5176,-1.9318){$+$};
\node [above]at(0.5176,1.9318){$+$};
\node [right]at(1.4142,-1.4142){$-$};
\node [above]at(-0.5176,1.9318){$+$};
\node [left]at(-1.9318,-0.5176){$-$};
\node [left]at(-1.4142,1.4142){$-$};
\node [left]at(-1.4142,-1.4142){$+$};
\node [left]at(-1.9318,0.5176){$+$};
\node [below]at(-0.5176,-1.9318){$-$};
\node [right]at(1.2,0.9){$a_1$};
\node [left]at(0.8,1){$b_1$};
\node [right]at(-1,1){$c_1$};
\node [right]at(-1,-1){$e_1$};
\node [right]at(-1.5,1.4){$d_1$};
\end{tikzpicture}}
\qquad
\qquad
\qquad
\subfigure[$G(D_2)$]{
  \begin{tikzpicture}[scale=0.8]
\draw[dashed](0,0) circle [radius=2];
\draw[thick](1.9318,0.5176) arc (15:345:2);
\draw[thick,->](1.4142,1.4142)--(-1.4142,1.4142);
\draw[thick,<-](0.5176,1.9318)--(-0.5176,-1.9318);
\draw[thick,->](-0.5176,1.9318)--(-1.4142,-1.4142);
\draw[thick,<-](-1.9318,0.5176)--(1.4142,-1.4142);
\draw[thick,->](-1.9318,-0.5176)--(0.5176,-1.9318);
\fill (1.9318,0.5176) circle [radius=0.075];
\fill (1.9318,-0.5176) circle [radius=0.075];
\node [right]at(1.4142,1.4142){$-$};
\node [above]at(0.5176,1.9318){$-$};
\node [above]at(-0.5176,1.9318){$+$};
\node [left]at(-1.4142,1.4142){$+$};
\node [left]at(-1.9318,0.5176){$+$};
\node [left]at(-1.9318,-0.5176){$-$};
\node [left]at(-1.4142,-1.4142){$-$};
\node [below]at(-0.5176,-1.9318){$+$};
\node [below]at(0.5176,-1.9318){$+$};
\node [right]at(1.4142,-1.4142){$-$};
\node [left]at(1.2,1.6){$a_2$};
\node [right]at(0.2,1){$b_2$};
\node [right]at(-0.75,1){$c_2$};
\node [right]at(-1,-1){$e_2$};
\node [left]at(1.2,-1.4){$d_2$};

\end{tikzpicture}}
  \caption{Gauss diagrams of $D_1$ and $D_2$.}
  \label{gd1}
\end{figure}
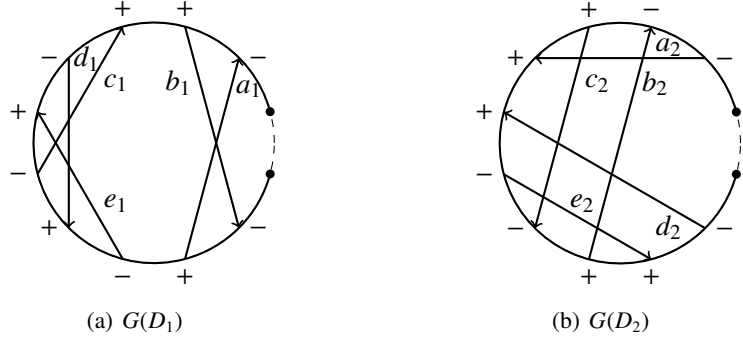

Firstly, the Gauss diagrams $G(D_1)$ and $G(D_2)$ of $D_1$  and $D_2$ are shown in Figure~\ref{gd1}.

For $D_1$ and $G(D_1)$, we can get ${\rm sign}(a_1)=1$, ${\rm sign}(b_1)=1$, ${\rm sign}(c_1)=-1$, ${\rm sign}(d_1)=-1$ and ${\rm sign}(e_1)=-1$.
$i(a_1)=-1$, $i(b_1)=-1$, $i(c_1)=0$, $i(d_1)=0$ and $i(e_1)=0$.

Then $F_{D_1}(t)=-2+2t^{-1}$.

 Next, we calculate the $n$th polynomial for $D_1$. For each $n\ge2$, $C_n(G(D_1))$= $\left\{
 c_1,~d_1,~e_1
 \right\}$ and $d_n(c_1)=d_n(d_1)=d_n(e_1)=0$. Then
 $$
  Z^n_{D_1}(t)=\left\{
  \begin{array}{cl}
  -2+2t^{-1}, & \ n=1,\\
  0, & \ n\ge2. \\
  \end{array}\right.
$$

For $D_2$ and $G(D_2)$, we can get ${\rm sign}(a_2)=-1$, ${\rm sign}(b_2)=1$, ${\rm sign}(c_2)=1$, ${\rm sign}(d_2)=-1$ and ${\rm sign}(e_2)=-1$. $i(a_2)=0$, $i(b_2)=-1$, $i(c_2)=-1$, $i(d_2)=0$ and $i(e_2)=0$.

Then $F_{D_2}(t)=-2+2t^{-1}$.

Next, we calculate the $n$th polynomial for $D_2$. For each $n\ge2$, $C_n(G(D_2))$= $\left\{
 a_2,~d_2,~e_2
 \right\}$ and $d_n(a_2)=d_n(d_2)=d_n(e_2)=0$. Then $$
  Z^n_{D_2}(t)=\left\{
  \begin{array}{cl}
  -2+2t^{-1}, &\  n=1,\\
  0, &\  n\ge2. \\
  \end{array}\right.
$$

 Finally, we calculate $F_{D_1}(u,v)$ and $F_{D_2}(u,v)$. According to the definition, we consider the subsets of chords that intersect with each of these five chords $x$ and the maps ${\phi_x}$ ($x=a_1,b_1,c_1,d_1,e_1$) of $D_1$. Then we have  $r(a_1)=\emptyset,\ l(a_1)=\left\{b_1\right\} \ {\rm and}\  \phi_{a_1}=0;\ r(b_1)=\emptyset, \ l(b_1)=\left\{a_1 \right\} \ {\rm and}\  \phi_{b_1}=0;\
 r(c_1)=\left\{
e_1
 \right\},\  l(c_1)=\left\{
 d_1
 \right\}\ {\rm and}\ \phi_{c_1}={\rm id};
 \ r(d_1)=\left\{
e_1
 \right\}, \ l(d_1)=\left\{
 c_1
 \right\} \ {\rm and}\ \phi_{d_1}={\rm id};\
 r(e_1)=\left\{
d_1
 \right\},\  l(e_1)=\left\{
 c_1
 \right\} \ {\rm and}\  \phi_{e_1}={\rm id}.$

 Next, we compute the index functions of the five crossings $x$ of $D_1$.
$$
    \begin{aligned}
  &g_{a_1}(v)=-{\rm sign}(b_1)v^{\phi_{a_1}(-i(b_1))}=-1,\\&
  g_{b_1}(v)=-{\rm sign}(a_1)v^{\phi_{b_1}(-i(a_1))}=-1,\\&
   g_{c_1}(v)={\rm sign}(e_1)v^{\phi_{c_1}(i(e_1))}-{\rm sign}(d_1)v^{\phi_{c_1}(-i(d_1))}=0,\\&
   g_{d_1}(v)={\rm sign}(e_1)v^{\phi_{d_1}(i(e_1))}-{\rm sign}(c_1)v^{\phi_{d_1}(-i(c_1))}=0,\\&
   g_{e_1}(v)={\rm sign}(d_1)v^{\phi_{e_1}(i(d_1))}-{\rm sign}(c_1)v^{\phi_{e_1}(-i(c_1))}=0.
   \end{aligned}
$$
 Then $F_{D_1}(u,v)=-2+2u^{-1}$.

Similarly, we consider the subsets of chords that intersect with each of these five chords $y$ and the maps ${\phi_y}$ ($y=a_2,b_2,c_2,d_2,e_2$) of $D_2$. Then we have
 $r(a_2)=\left\{
b_2
 \right\}, \ l(a_2)=\left\{
c_2
 \right\}\  {\rm and}\  \phi_{a_2}={\rm id};
\ r(b_2)=\left\{
a_2, d_2
 \right\}, \ l(b_2)=\left\{
 e_2
 \right\}\  {\rm and}\   \phi_{b_2}=0;
 \ r(c_2)=\left\{
a_2,d_2
 \right\}, \ l(c_2)=\left\{
e_2
 \right\} \  {\rm and}\   \phi_{c_2}=0;
  \ r(d_2)=\left\{
c_2
 \right\},\  l(d_2)=\left\{
 b_2
 \right\} \  {\rm and}\   \phi_{d_2}={\rm id};
  \ r(e_2)=\left\{
c_2
 \right\},\  l(e_2)=\left\{
 b_2
 \right\} \  {\rm and}\   \phi_{e_2}={\rm id}.$

Next, we compute the index functions of the five crossings $y$ of $D_2$.
$$
    \begin{aligned}
  &g_{a_2}(v)={\rm sign}(b_2)v^{\phi_{a_2}(i(b_2))}-{\rm sign}(c_2)v^{\phi_{a_2}(-i(c_2))}=v^{-1}-v,\\&
  g_{b_2}(v)={\rm sign}(a_2)v^{\phi_{b_2}(i(a_2))}+{\rm sign}(d_2)v^{\phi_{b_2}(i(d_2))}-{\rm sign}(e_2)v^{\phi_{b_2}(-i(e_2))}=-1,\\&
   g_{c_2}(v)={\rm sign}(a_2)v^{\phi_{c_2}(i(a_2))}+{\rm sign}(d_2)v^{\phi_{c_2}(i(d_2))}-{\rm sign}(e_2)v^{\phi_{c_2}(-i(e_2))}=-1,\\&
   g_{d_2}(v)={\rm sign}(c_2)v^{\phi_{d_2}(i(c_2))}-{\rm sign}(b_2)v^{\phi_{d_2}(-i(b_2))}=v^{-1}-v,\\&
   g_{e_2}(v)={\rm sign}(c_2)v^{\phi_{e_2}(i(c_2))}-{\rm sign}(b_2)v^{\phi_{e_2}(-i(b_2))}=v^{-1}-v.
   \end{aligned}
$$
  Then $F_{D_2}(u,v)=-3u^{v^{-1}-v}+2u^{-1}+1$.

To sum up, $F_{D_1}(t)=F_{D_2}(t)$ and $Z^n_{D_1}(t)=Z^n_{D_2}(t)$, but $F_{D_1}(u,v)\neq F_{D_2}(u,v)$. So $F$-polynomial can detect $D_1$ is not equivalent to $D_2$ while the index polynomial and $n$th  polynomial cannot.

\begin{figure}
 \centering
\begin{tikzpicture}[scale=0.7]
\draw[thick,->](0,0)--(0,0.5);
\draw[thick](0,0.5)--(0,1.9);
\draw[thick](0,2.1)--(0,3);
\draw[thick](0,3)--(-2,3);
\draw[thick](-2,3)--(-2,1);
\draw[thick](-2,1)--(-0.1,1);
\draw[thick](0.1,1)--(8,1);
\draw[thick](8,1)--(8,4);
\draw[thick](8,4)--(7.1,4);
\draw[thick](6.9,4)--(5.1,4);
\draw[thick](2.5,4)--(2.1,4);
\draw[thick](4.9,4)--(3.5,4);
\node[right]at(2.5,4){$\cdots$};
\draw[thick](1.9,4)--(-1,4);
\draw[thick](-1,4)--(-1,3.1);
\draw[thick](-1,2.9)--(-1,2);
\draw[thick](-1,2)--(1,2);
\draw[thick](1,2)--(1,3.9);
\draw[thick](1,4.1)--(1,5);
\draw[thick](1,5)--(2,5);
\draw[thick](2,5)--(2,3);
\draw[thick](2,3)--(2.5,3);
\node[above]at(2.5,3){$\vdots$};
\draw[thick](4,3)--(3.5,3);
\node[above]at(3.5,3){$\vdots$};
\draw[thick](7,3)--(7,5);
\draw[thick](7,5)--(6,5);
\draw[thick](6,5)--(6,4.1);
\draw[thick](6,3.9)--(6,3);
\draw[thick](6,3)--(5,3);
\draw[thick](5,3)--(5,5);
\draw[thick](5,5)--(4,5);
\draw[thick](4,5)--(4,4.1);
\draw[thick](4,3.9)--(4,3);
\fill (0,0) circle (.1);
\fill (7,3) circle (.1);
\node [below] at (-0.2,1){$a$};
\node [below] at (-0.2,2){$b$};
\node [above] at (-1.2,3){$c$};
\node [above] at (7.3,4){$d_1$};
\node [above] at (6.3,4){$e_1$};
\node [above] at (5.3,4){$d_2$};
\node [above] at (4.3,4){$e_2$};
\node [above] at (1.35,4){$e_m$};
\node [above] at (2.35,4){$d_m$};
\end{tikzpicture}
\caption{A family of knotoid diagrams $H_m$.}
\label{a family}
\end{figure}
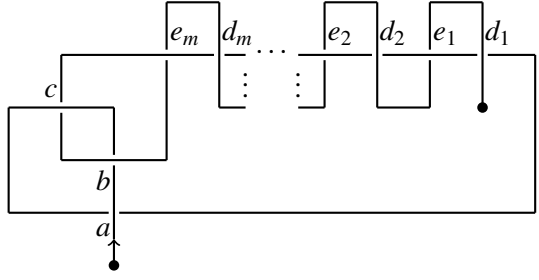

\begin{example}\label{hm}
{\rm Let $\left\{H_m\right\}$ be a family of knotoid diagrams and $m\ge2$ as shown in Figure~\ref{a family}.}
\end{example}

\begin{figure}
 \centering
\begin{tikzpicture}[scale=0.9]
\draw[dashed](0,0) circle [radius=2];
 \draw[thick] (1.732,1) arc (30:375:2);
 \fill (1.9318,0.5176) circle [radius=0.075];
\fill (1.732,1) circle[radius=0.075];
\draw[thick,->](1.4142,1.4142)--(-1.4142,1.4142);
\draw[thick,->](0.5176,-1.9318)--(0.5176,1.9318);
\draw[thick,->](-0.5176,1.9318)--(-0.5176,-1.9318);
\draw[thick,->](2,0)--(-2,0);
\draw[thick,->](-1.9318,-0.5176)--(1.9318,-0.5176);
\draw[thick,->](-1,-1.732)--(1,-1.732);
\draw[thick,->](1.4142,-1.4142)--(-1.4142,-1.4142);
\node [above]at(0,1.4){$a$};
\node [right]at(0.5176,1){$b$};
\node [left] at(-0.5176,1){$c$};
\node [above]at(0,-0.05){$d_1$};
\node [above]at(0,-1.85){$e_m$};
\node [above]at(0,-0.55){$e_1$};
\node [above]at(0,-1.4142){$d_m$};
\node [right]at(0.5176,-0.9){\vdots};
\node [left]at(-1.4142,1.4142){$+$};
\node [right]at(1.4142,1.4142){$-$};
\node [above]at(0.5176,1.9318){$-$};
\node [below]at(0.5176,-1.9318){$+$};
\node [below]at(-0.5176,-1.9318){$-$};
\node [above]at(-0.5176,1.9318){$+$};
\node [left]at(-2,0){$+$};
\node [right]at(2,0){$-$};
\node [right]at(1.9318,-0.5176){$+$};
\node [left]at(-1.9318,-0.5176){$-$};
\node [left]at(-1.1,-1.732){$-$};
\node [right]at(1,-1.75){$+$};
\node [right]at(1.4142,-1.4142){$-$};
\node [left]at(-1.4142,-1.4142){$+$};
\end{tikzpicture}
\caption{Gauss diagrams $G(H_m)$ of $H_m$.}
\label{a family g}
\end{figure}
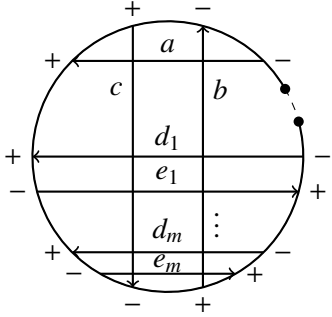

Firstly, we give the Gauss diagrams of $H_m$ as shown in Figure~\ref{a family g}.

We have ${\rm sign}(a)=-1$, ${\rm sign}(b)=1$, ${\rm sign}(c)=1$, ${\rm sign}(d_1)={\rm sign}(d_2)={\cdots}={\rm sign}(d_m)=-1$ and ${\rm sign}(e_1)={\rm sign}(e_2)={\cdots}={\rm sign}(e_m)=-1$. $i(a)=0$, $i(b)=-1$, $i(c)=-1$, $i(d_1)=i(d_2)={\cdots}=i(d_m)=0$ and $i(e_1)=i(e_2)={\cdots}=i(e_m)=0$.

Then the index polynomial $F_{H_m}(t)=-2+2t^{-1}$.

Next, we calculate the $n$th polynomial of $H_m$. For each $n\ge2$, $C_n(G(H_m))= \left\{
 a,~d_1,~d_2,{\cdots},~d_m,~e_1,~e_2,\right.$
 $\left.{\cdots},~e_m
 \right\}$, $d_n(a)=d_n(d_1)={\cdots}=d_n(d_m)=d_n(e_1)={\cdots}=d_n(e_m)=0$. Then
$$
   Z^n_{H_m}(t)=\left\{
  \begin{array}{cl}
  -2+2t^{-1}, & \ n=1,\\
  0, &\  n\ge2 .\\
  \end{array}\right.
$$
  But we have $g_{a}(v)=v^{-1}-v,\ g_{b}(v)=-1,\
   g_{c}(v)=-1,\
   g_{d_1}(v)=g_{d_2}(v)={\cdots}=g_{d_m}(v)=v^{-1}-v ~{\rm and}~
  g_{e_1}(v)=g_{e_2}(v)={\cdots}=g_{e_m}(v)=v^{-1}-v.$

  Then we have
  $F_{H_m}(u,v)=-(2m+1)u^{v^{-1}-v}+2u^{-1}+2m-1.$

 For different integers $m_1$ and $m_2$ that are greater than or equal to 2, we have $2m_1\pm1 \neq 2m_2\pm1$, so $ F_{H_{m_1}}(u,v) \neq F_{H_{m_2}}(u,v)$. It follows that
   $H_m$ are distinct from each other.

According to Example~\ref{hm}, we do construct a family of knotoid diagrams that the index polynomial and the $n$th polynomial cannot distinguish while the $F$-polynomial can.

\end{document}